\documentclass[12pt,reqno]{amsart}
\usepackage{amssymb,latexsym}
\usepackage{amsmath}
\usepackage{amsthm}
\usepackage{graphicx}
\usepackage{hyperref}
\usepackage{titletoc}
\usepackage{pdfsync}
\usepackage{color,xcolor}
\usepackage{amsthm,amsmath,amssymb}

\usepackage{mathrsfs}

\numberwithin{equation}{section}
\newtheorem{theorem}{Theorem}[section]
\newtheorem{proposition}[theorem]{Proposition}
\newtheorem{lemma}[theorem]{Lemma}

\usepackage{multirow}
\usepackage[font=bf,aboveskip=15pt]{caption}
\usepackage[toc,page]{appendix}

\theoremstyle{definition}

\newtheorem{remark}[theorem]{Remark}

\def\C{{\mathfrak C}}

\def\R{{\mathfrak R}}

\def\begeq{\begin{equation}}
	\def\endeq{\end{equation}}

\def\R{\Bbb R}

\allowdisplaybreaks

\begin{document}

\title[New type of solutions for Schr\"odinger equations with critical growth
 ]
{ New type of solutions for Schr\"odinger equations with  critical growth
 }

 \author{Yuan Gao and  Yuxia Guo
}

\address{Yuan Gao,
		\newline\indent Department of Mathematical Sciences, Tsinghua University,
		\newline\indent Beijing 100084,  P. R. China.
	}
	\email{gaoy22@mails.tsinghua.edu.cn}

	\address{Yuxia Guo,
		\newline\indent Department of Mathematical Sciences, Tsinghua University,
		\newline\indent Beijing 100084,  P. R. China.
	}
	\email{yguo@tsinghua.edu.cn}

\begin{abstract}

We  consider the following  nonlinear Schr\"odinger equations with critical growth:
\begin{equation}
\label{eq}
- \Delta u + V(|y|)u=u^{\frac{N+2}{N-2}},\quad u>0 \ \   \mbox{in} \  \R^N, \ \ \
\end{equation}
where $V(|y|)$ is a bounded positive radial function in $C^1$, $N\ge 5$. By using a finite reduction argument, we
show that if $r^2V(r)$ has either an isolated local maximum or an isolated minimum at $r_0>0$ with $V(r_0)>0$, there exists infinitely many non-radial large energy solutions which are invariant under some sub-groups of $O(3)$.
\vspace{2mm}
		
		
		\vspace{2mm}
		

	\end{abstract}

\maketitle
\section{Introduction}\label{sec1}
We consider the standing wave solutions of the form
$
\psi(y,t) = e^{{  \bf{ i } } \lambda t } u(y)
$ for the following time-dependent Schr\"odinger equation
	\begin{align}\label{original1}
		-{  \bf{ i } } \frac{\partial \psi}{\partial t}  = \Delta \psi - \widetilde V(y) \psi + |\psi|^{p-1} \psi,     \quad\text{for} ~ (y,t) \in \mathbb{R}^N\times \mathbb{R},
	\end{align}
where ${\bf{i}}$ is the imagine,  $p>1.$ 
Under the assumption that $u(y)$ is positive and vanishes at infinity, we see that $\psi(y,t)$  is a solution of \eqref{original1} if and only if $u(y)$ is the solution of  the following nonlinear elliptic equation
\begin{equation}\label{equation0}
-\Delta u + V(y) u = u^p, \quad u>0 \quad\text{in}\; \mathbb R^N,\quad
\lim\limits_{|y|\to+\infty} u(y)=0,
\end{equation}
where $\widetilde V(y)=V(y)-\lambda.$
\\

For the subcritical case, that is  $p\in \left(1,\displaystyle\frac{N+2}{N-2}\right),$ there have been a lot of existence results obtained in the literature, see for example, \cite{dWY,DM,GG,GMPY1,WY1} and the references therein. 
For the critical case, that is $p=\displaystyle\frac{N+2}{N-2}$, \eqref{equation0} turns to be:
\begin{equation}\label{equation}
\begin{cases}
-\Delta u + V(y) u = u^{\frac{N+2}{N-2}}, \quad u>0 & \text{in}\; \mathbb R^N,\\[2mm]
\lim\limits_{|y|\to+\infty} u(y)=0,
\end{cases}
\end{equation}
with the assumptions $V(y)\ge 0$ and $V(y) \not\equiv 0.$ Problem \eqref{equation} is also closely related to the well known Brezis-Nirenberg problem on $\mathbb{S}^N:$
\begin{equation}\label{BN}
    -\Delta_{\mathbb{S}^N}u=u^{\frac{N+2}{N-2}}+\lambda u,\quad u>0\quad \text{on}\quad \mathbb{S}^N,
\end{equation}
which if $\lambda<-\displaystyle\frac{N(N-2)}{4}$, after the stereo-graphic projection, can be reduced to \eqref{equation} with
\begin{equation*}
    V(y)=\frac{-4\lambda-N(N-2)}{(1+|y|^2)^2}.
\end{equation*}
For more details about \eqref{BN}, we refer readers to \cite{BW,BL,BP,D,D1,D2,GS} and the references therein.
In particular, we see that if $V(y)\ge 0$ and $V(y) \not\equiv 0,$ the mountain pass value for \eqref{equation} is not a critical value for the corresponding functional. Thus we can not use concentration compactness arguments (see \cite{L,L1}) to obtain a solution for \eqref{equation}. In \cite{WY}, Chen, Wei and Yan firstly constructed infinitely many non-radial positive solutions for \eqref{equation} in the radially symmetric potential case, that is $V(|y|)\ge 0$ radially symmetric and bounded, with the additional assumption that $r^2V(r)$ has a local maximum or minimum point $r_0>0, V(r_0)>0$. Their solutions are obtained by gluing a large number of Aubin-Talenti bubbles (see \cite{T}):
\begin{equation}
    U_{x,\mu}(y)=C_N\left(\frac{\mu}{1+\mu^2|y-x|^2}\right)^{\frac{N-2}{2}},\quad C_N=(N(N-2))^{\frac{N-2}{4}},
\end{equation}
at the circle of $(x_1,x_2)$-plane:
\begin{equation}
    x_j^*=r\Bigl(\cos\frac{2(j-1)\pi}k, \sin\frac{2(j-1)\pi}k,{\bf 0}\Bigr),\quad j=1,\cdots,k,
\end{equation}
where ${\bf 0}$ is the zero vector in $\R^{N-2}, k$ is a large integer, $r \sim r_0, \mu \sim k^{\frac{N-2}{N-4}}$ as $k \to +\infty.$ 
And the main order of the solutions built in \cite{WY} looks like
\begin{equation*}
    \bar u_k \sim \sum_{j=1}^k U_{x_j^*,\bar \mu}.
\end{equation*}
This kind of solutions is usually called bubble solutions or equatorial solutions (see \cite{DM}). After that, Peng, Wang and Yan \cite{PWY} tried to weaken the radially symmetry condition for $V(y)$ to $V(y)=V(|y'|,y''),$ where $(y',y'')\in\R^2\times\R^{N-2},$ and successfully constructed the bubble solutions concentrated at the critical point of $r^2V(r,y'')$.
Recently, He, Wang and Wang \cite{HWW} proved a non-degeneracy result for the solution found in \cite{PWY}. And as an application of non-degeneracy result, they constructed a new type of solutions by gluing large number of bubbles located at the circle of $(x_1, x_2) $-plane and $(x_3, x_4) $-plane separately. For any other related work about bubble solutions, we refer the reader to \cite{GMPY,LWX,MWY,WY2,WY3} and the references therein.

In this paper, we will construct a totally different type of solutions to \eqref{equation} with more complex concentration structure, which we called, double-tower solutions. This kind of solutions are centered at the points lying on the top and the bottom circles of a cylinder with height $h$. More precisely, set $k$ be an integer number,
	\begin{equation*}
		x_{j}^{+}=r\Bigl(\sqrt{1-h^2}\cos\frac{2(j-1)\pi}k, \sqrt{1-h^2}\sin\frac{2(j-1)\pi}k,h,{\bf 0}\Bigr),\quad j=1,\cdots,k,
	\end{equation*}
	and
	\begin{equation*}
		x_{j}^{-}=r\Bigl(\sqrt{1-h^2}\cos\frac{2(j-1)\pi}k, \sqrt{1-h^2}\sin\frac{2(j-1)\pi}k,-h,{\bf 0}\Bigr),\quad j=1,\cdots,k,
	\end{equation*}
	where ${\bf 0}$ is the zero vector in $\R^{N-3}$. Then the solutions we built seem like:
 \begin{equation}
     u_k \sim W_{r,h,\mu}(y):=\sum_{j=1}^k U_{x_j^+,\mu}(y)+ \sum_{j=1}^k U_{x_j^-,\mu}(y).
 \end{equation}
 Indeed, there have been several existence results for such kind of double-tower solutions to other related equations, for example, \cite{DMW,GGH,MM} and the references therein.

To present our results, we assume the potential function $V(y)$ satisfies the following condition:
\medskip

$({\bf V}): V(y)=V(|y|) \ge 0 $ is  radially symmetric, positive, bounded and belongs to $C^1$. Moreover, $r^2V(r)$ has either an isolated local maximum or an isolated minimum at $r_0>0$ with $V(r_0)>0$.

For any point $y \in \mathbb{R}^{N}$, we set $y = (y',y''),y'\in\mathbb{R}^{3},y''\in\mathbb{R}^{N-3}$. Define
	\[
	\begin{split}
		H_{s}=\bigg\{ u: \ &u\in \overline{H^1(\R^N) \cap D^{1,2}(\R^N)} , \quad u \ \text{is even in} \ y_2,y_4,y_5,\cdots, y_N, \\ u&\left(\sqrt{y_1^2+y_2^2}\cos\theta , \sqrt{y_1^2+y_2^2}\sin\theta,y_3, y''\right)
		 \\&=u\left(\sqrt{y_1^2+y_2^2}\cos\left(\theta+\frac{2\pi j}k\right) , \sqrt{y_1^2+y_2^2}\sin\left(\theta+\frac{2\pi j}k\right), y_3,y''\right)
		\bigg\}.
	\end{split}
	\]
 Throughout of the this paper, we suppose that $N\ge5$ and $(r,h,\mu)\in\mathscr{S}_k,$ where
 \begin{align}\label{Sk}
			{{\mathscr S}_k}
			= \Bigg\{&(r,h,\mu) :\,  r\in \Big[r_0-\hat \sigma, r_0+\hat  \sigma \Big], \quad
			h \in \Bigg[\frac{h_0-\hat \sigma}{k^{\frac{N-3}{N-1}}}  ,
			\frac{h_0+\hat \sigma}{k^{\frac{N-3}{N-1}}} \Bigg], \nonumber
			\\[2mm]
			& \qquad \qquad
			\mu \in \Big[(\mu_0-\hat\sigma) k^{\frac{N-2}{N-4}},  (\mu_0+\hat \sigma) k^{\frac{N-2}{N-4}}\Big]\Bigg\},
		\end{align}
		with  $ \hat \sigma$ is a small fixed number, independent of $k$. 

\vskip8pt

Our main results are the following:
  \begin{theorem}\label{main1}
      Suppose that $V(y)$  satisfies $({\bf V})$ and $N\ge 5$. Then there is an
		integer $k_0>0$, such that for any integer  $k\ge k_0$, \eqref{equation}
		has a solution $u_k$ of the form
		\[
		u_k = W_{r_k,h_k,\mu_k}(y)+\phi_k,
		\]
		where  $(r_k,h_k,\mu_k) \in \mathscr{S}_k,\phi_k\in H_s$,
		and as $k\to +\infty$,
		$
		 ||\phi_k ||_{D^{1,2}(\R^N)}\to 0.
		$
  \end{theorem}
  
  	As a consequence, we have

\begin{theorem}\label{main2}
		Under the same assumptions as in Theorem \ref{main1},  \eqref{equation} has infinitely many positive non-radial solutions, which  are invariant under some non-trivial sub-groups of $O(3),$ and their energy can be made arbitrarily large.
	\end{theorem}
 
	Let us briefly outline the main idea for the proof of Theorem \ref{main1}. With the approximation solution $W_{r,h,\mu}$ in mind. Our main task is to find the correction term $\phi$ so that $u=W_{r,h,\mu}+\phi$ solves \eqref{equation}. Recall that \cite{R} the following functions
 \begin{equation}\label{Zj}
	Z_i(y) :={\partial U \over \partial y_i}(y), \quad i=1, \ldots , N, \quad
	Z_{N+1}(y) :={N-2 \over 2} U(y)+y\cdot \nabla U(y),
\end{equation}
span the kernel of the linear operator associated to the limit equation 
$$-\Delta u=u^{\frac{N+2}{N-2}} \hbox { in }\mathbb{R}^N,$$
that is the set of the solution to
\begin{equation}
    -\Delta \phi =(2^*-1)U^{2^*-2}\phi,\quad \phi\in D^{1,2}(\R^N).
\end{equation}
We first consider the linearlized problem \eqref{lin} about $\phi,$ so that problem \eqref{equation} can be reduced to a finite dimensional problem. Solve this auxiliary problem, we can obtain $\phi=\phi_{r,h,\mu}$. Next, we consider the critical point of the energy functional corresponding to \eqref{equation}:
\begin{equation}
    I(u)=\frac{1}{2}\int_{\R^N}\left(|Du|^2+V(|y|)u^2\right)-\frac{1}{2^*}\int_{\R^N}|u|^{2^*},
\end{equation}
to choose the proper parameters $(r,h,\mu)$ so that we can get a true solution to \eqref{equation}. Compared with the construction of bubble solutions, we need to overcome the difficulty caused by the third parameter $h.$ In fact, expansions for the main term of energy functional are determined by the order of $h$. Thus we should first decide the order of $h$. Roughly speaking, we have three cases:\\
  $\bullet${ \bf Case  1.} $h=o_k(1);$
			\\[2mm]
 $\bullet${ \bf Case  2.}  $\sqrt{1-h^2}=o_k(1);$
			\\[2mm]
$\bullet$ {\bf Case  3.}   $h\neq o_k(1),\sqrt{1-h^2}\neq o_k(1);$\\[2mm]
 However, in {\bf Case~2}, we can not find a critical point to the main term of the energy functional, and {\bf Case~3} is quite similar to the bubble solution case in \cite{WY}. So we study on  {\bf Case~1} and finally determine that $h\sim\displaystyle\frac{1}{k^{\frac{N-3}{N-1}}}.$

This paper is organized as follows. In Section \ref{secA}, we will calculate the expansions for the energy functional  $I(W_{r,h,\mu})$ and its partial derivatives to discuss the order of parameters which allows the existence of a critical point to $I(W_{r,h,\mu})$. In Section \ref{sec2}, we establish the linearized theory for the linearized problem and give estimates for the correction terms. Section \ref{sec3} is devoted to the proof of Theorem \eqref{equation}.  We will show the existence of a critical point to $I(W_{r,h,\mu}+\phi_{r,h,\mu})$ with some proper parameters $(r,h,\mu)\in \mathscr{S}_k.$  We put some basic estimates and useful lemmas  in Appendix \ref{secB}.

For the reader's convenience, we will provide a collection of notations. Throughout this paper, we employ  $C $ to denote certain constants and  $ \sigma, \hat \sigma, \tilde\sigma, \delta, \tau, \eta, \eta_0  $ to denote some small constants or functions.  We also use $ \delta_{ij} $  to denote Kronecker delta function:
	\[ \delta_{ij}= \begin{cases} 1,  \quad \text{if} ~ i= j, \\[2mm]
		0,  \quad \text{if} ~ i \neq j.
	\end{cases} \]
	Furthermore, we  also  employ the notation  by writing  $O(f(r,h,\mu)), o(f(r,h,\mu))  $ for the functions which satisfy
	\begin{equation*}
		\text{if} \quad g(r,h,\mu) \in O(f(r,h,\mu)) \quad \text{then}\quad {\lim_{k \to+\infty}} \Bigg|\, \frac{g(r,h,\mu)}{f(r,h,\mu)} \, \Bigg|\leq C<+\infty,
	\end{equation*}
	and
	\begin{equation*}
		\text{if} \quad g(r,h,\mu) \in o(f(r,h,\mu)) \quad \text{then}\quad {\lim_{k \to+\infty}} \frac{g(r,h,\mu)}{f(r,h,\mu)}=0.
	\end{equation*}
\medskip

\section{Expansions for the Energy functional $I(W_{r,h,\mu})$}\label{secA}
This section is devoted to the energy functional  $I(W_{r,h,\mu})$ and its partial derivatives.  We first give the following Lemma.
\begin{lemma}\label{A.3}
Assume that $V(|y|)$ satisfies $({\bf V})$ and $N\ge 5, \displaystyle\frac{k}{\mu\sqrt{1-h^2}}=o(1), \mu h \to +\infty$, we have the following expansion of the energy functional for $k \to +\infty$:
\begin{equation}\label{I(W)}
\begin{aligned}
    I(W_{r,h,\mu})=&k\left(A_1+\frac{A_2 V(r)}{\mu^2}-\left(\sum_{j=2}^{k}\frac{B_0}{\mu^{N-2}|{x}_{j}^{+}-{x}_{1}^{+}|^{N-2}}+\sum_{j=1}^{k}\frac{B_0}{\mu^{N-2}|{x}_{1}^{+}-{x}_{j}^{-}|^{N-2}}\right)\right)\\
    &+ko\left(\frac{1}{\mu^2}\right)+k O\left(\frac{1}{\mu^N d_0^N}\right)+kO\left( \frac{1}{\mu^{N-2}}\left(\sum_{j=2}^{k}\frac 1 {|{x}_{j}^{+}-{x}_{1}^{+}|^{N-4}}+\sum_{j=1}^{k}\frac 1 {|{x}_{1}^{+}-{x}_{j}^{-}|^{N-4}}\right)  \right)\\
        &+\frac{k}{\mu^{N-\sigma}}O\left(\sum_{j=2}^{k}\frac 1 {|{x}_{j}^{+}-{x}_{1}^{+}|^{N-\sigma}}+\sum_{j=1}^{k}\frac 1 {|{x}_{1}^{+}-{x}_{j}^{-}|^{N-\sigma}}\right)\\
&+\frac{k}{\mu^{N-\sigma}}O \left(\sum_{j=2}^{k}\frac{1}{|{x}_{j}^{+}-{x}_{1}^{+}|^\frac{N-\sigma}{2}}+\sum_{j=1}^{k}\frac{1}{|{x}_{1}^{+}-{x}_{j}^{-}|^\frac{N-\sigma}{2}}\right)^2
\\&+\frac{k}{\mu^{N-\sigma}}O\left(\sum_{j=2}^{k}\frac{1}{|{x}_{j}^{+}-{x}_{1}^{+}|^\frac{(N-\sigma)(N-2)}{2N}}+\sum_{j=1}^{k}\frac{1}{|{x}_{1}^{+}-{x}_{j}^{-}|^\frac{(N-\sigma)(N-2)}{2N}}\right)^{2^*},
\end{aligned}
\end{equation}
where $A_1=\displaystyle\frac{2}{N} \int_{\R^N}U_{0, 1}^{2^*}, A_2=\int_{\R^N}U_{0, 1}^{2}, B_0=\displaystyle\int_{\mathbb{R}^N}  \frac{C_N^{2^*}}{(1+z^{2})^{\frac{N+2} 2}}=\frac{N}{2}A_1 $ are positive constants and  $ \sigma $ is a  constant small enough.
\end{lemma}

\begin{proof}
For  $j= 1,\cdots, k $, let
	\begin{align*}
		\Omega_j := &\Big\{y=(y_1, y_2, y_3, y'') \in \mathbb{R}^3 \times \mathbb{R}^{N-3}:
		\nonumber\\[2mm]
		& \qquad\Big\langle \frac{(y_1, y_2)}{|(y_1, y_2)|},  \Big(\cos{\frac{2(j-1) \pi}{k}}, \sin{\frac{2(j-1) \pi}{k}}\Big)  \Big\rangle\geq \cos{\frac \pi k}\Big\}.
	\end{align*}
 Furthermore, we divide $\Omega_j$ into two parts:
	\begin{align*}
		\Omega_j^+= & \Big\{y:  y=(y_1, y_2, y_3, y'')  \in  \Omega_j, y_3\geq0 \Big\},
		\\[2mm]
		\Omega_j^-= & \Big\{y: y=(y_1, y_2, y_3, y'')  \in  \Omega_j, y_3<0 \Big\},
	\end{align*}
	then
	$$\mathbb{R}^N= \underset{j=1}{\overset{k}{\cup}} \Omega_j,  \quad \Omega_j=  \Omega_j^+\cup  \Omega_j^-$$
	and
	$$ \Omega_j \cap  \Omega_i=\emptyset,  \quad  \Omega_j^+\cap  \Omega_j^-=\emptyset, \qquad \text{if}\quad i\neq j.$$

Recall that $$I(u)=\frac{1}{2}\int_{\R^N}\left(|Du|^2+V(|y|)u^2\right)-\frac{1}{2^*}\int_{\R^N}|u|^{2^*}.$$
Using the symmetry, we obtain
\begin{equation}\label{I}
    \begin{aligned}
&I(W_{r,h,\mu})\\=&\frac{1}{2}\int_{\R^N}\left(|DW_{r,h,\mu}|^2+V(|y|)W_{r,h,\mu}^2\right)-\frac{1}{2^*}\int_{\R^N}|W_{r,h,\mu}|^{2^*}\\
=&k\int_{\R^N} U_{0,1}^{2^*}+k\left(\sum_{j=2}^{k}\int_{\R^N}  U_{{x}_{1}^{+}, \mu}^{2^*-1}  U_{{x}_{j}^{+}, \mu}+\sum_{j=1}^{k}\int_{\R^N}  U_{{x}_{1}^{+}, \mu}^{2^*-1}  U_{{x}_{j}^{-}, \mu}\right)+k\int_{\R^N} V(|y|)U_{{x}_{1}^{+}, \mu}^{2}\\
&+k\left(\sum_{j=2}^{k}\int_{\R^N} V(|y|)U_{{x}_{1}^{+}, \mu}U_{{x}_{j}^{+}, \mu}+\sum_{j=1}^{k}\int_{\R^N} V(|y|)U_{{x}_{1}^{+}, \mu}U_{{x}_{j}^{-}, \mu}\right)-\frac{1}{2^*}\int_{\R^N}|W_{r,h,\mu}|^{2^*}\\
:=&kI_0+kI_1+kI_2+kI_3-\frac{1}{2^*}I_4.
    \end{aligned}
\end{equation}
    Next, we calculate term by term. For $I_1$, by Lemma \ref{A.1} we have
    \begin{equation}\label{I1}
    \begin{aligned}
        I_1=&\left(\sum_{j=2}^{k}\int_{\R^N}  U_{{x}_{1}^{+}, \mu}^{2^*-1}  U_{{x}_{j}^{+}, \mu}+\sum_{j=1}^{k}\int_{\R^N}  U_{{x}_{1}^{+}, \mu}^{2^*-1}  U_{{x}_{j}^{-}, \mu}\right)\\
        =&\frac{B_0}{\mu^{N-2}}\left(\sum_{j=2}^{k}\frac{1}{|{x}_{j}^{+}-{x}_{1}^{+}|^{N-2}}+\sum_{j=1}^{k}\frac{1}{|{x}_{1}^{+}-{x}_{j}^{-}|^{N-2}}\right)\\
        &+O\left(\frac{1}{\mu^{N-\sigma}}\left(\sum_{j=2}^{k}\frac 1 {|{x}_{j}^{+}-{x}_{1}^{+}|^{N-\sigma}}+\sum_{j=1}^{k}\frac 1 {|{x}_{1}^{+}-{x}_{j}^{-}|^{N-\sigma}}\right)\right).
    \end{aligned}
    \end{equation}
Using Lemma \ref{A.2}, we also have
\begin{equation}\label{I3}
    I_3=O\left(\frac{1}{\mu^{N-2}}\left(\sum_{j=2}^{k}\frac 1 {|{x}_{j}^{+}-{x}_{1}^{+}|^{N-4}}+\sum_{j=1}^{k}\frac 1 {|{x}_{1}^{+}-{x}_{j}^{-}|^{N-4}}\right)\right).
\end{equation}
For $I_2$, for any small fixed $\delta>0$, we have
\begin{equation}\label{I2}
\begin{aligned}
    I_2=&\int_{\R^N} V(|y|)U_{{x}_{1}^{+}, \mu}^{2}=\int_{B_{\delta}({x}_{1}^{+})} V(|y|)U_{{x}_{1}^{+}, \mu}^{2}+\int_{\R^N \setminus B_{\delta}({x}_{1}^{+})} V(|y|)U_{{x}_{1}^{+}, \mu}^{2}\\=&\frac{1}{\mu^2}\int_{B_{\mu\delta}(0)}V(|\mu^{-1}z+{x}_{1}^{+}|)U_{0, 1}^{2}+O\left(\frac{1}{\mu^{N-2}}\right)=\frac{1}{\mu^2}\left(\int_{\R^N}V(r)U_{0, 1}^{2}+o(1)\right).
\end{aligned}
\end{equation}
Finally, we consider $I_4$:
\begin{equation}
    \begin{aligned}
I_4=&\int_{\R^N}|W_{r,h,\mu}|^{2^*}=2k\int_{\Omega_{1}^{+}}|W_{r,h,\mu}|^{2^*}\\
=&2k\left(\int_{\R^N}U_{0, 1}^{2^*}+\int_{\R^N \setminus \Omega_{1}^{+}}U_{{x}_{1}^{+}, \mu}^{2^*}+2^* I_1+\int_{\R^N \setminus \Omega_{1}^{+}}U_{{x}_{1}^{+}, \mu}^{2^*-1}\left(\sum_{j=2}^{k}U_{{x}_{j}^{+}, \mu}+\sum_{j=1}^{k}U_{{x}_{j}^{-}, \mu}\right)\right)\\
&+kO\left(\int_{\Omega_{1}^{+}}U_{{x}_{1}^{+}, \mu}^{2^*-2}\left(\sum_{j=2}^{k}U_{{x}_{j}^{+}, \mu}+\sum_{j=1}^{k}U_{{x}_{j}^{-}, \mu}\right)^2+\int_{\Omega_{1}^{+}}\left(\sum_{j=2}^{k}U_{{x}_{j}^{+}, \mu}+\sum_{j=1}^{k}U_{{x}_{j}^{-}, \mu}\right)^{2^*}\right).\\
    \end{aligned}
\end{equation}
We can also check that
\begin{equation*}
    \int_{\R^N \setminus \Omega_{1}^{+}}U_{{x}_{1}^{+}, \mu}^{2^*}=O\left(\int_{B_{d_0}^C({x}_{1}^{+})}\frac{\mu^N}{(1+\mu^2
|y-{x}_{1}^{+}|^2)^N}\right)=O\left(\frac{1}{\mu^Nd_0^N}\right),
\end{equation*}
and from \eqref{a12}--\eqref{a13}, we have
\begin{equation*}
    \int_{\R^N \setminus \Omega_{1}^{+}}U_{{x}_{1}^{+}, \mu}^{2^*-1}U_{{x}_{j}^{+}, \mu}=O\left(\frac{1}{\mu^Nd_0^N}\right),
\end{equation*}
where $d_0:=\min\left\{\displaystyle\frac{1}{2}|{x}_{2}^{+}-{x}_{1}^{+}|,\frac{1}{2}|{x}_{1}^{+}-{x}_{1}^{-}|  \right\}$. \\
On the other hand, by Lemma \ref{B.3},  we have that for $y \in \Omega_1^{+}$,
    \begin{equation}
        \begin{aligned}
\sum_{j=2}^{k}U_{{x}_{j}^{+},\mu}+\sum_{j=1}^{k}U_{{x}_{j}^{-},\mu} \le \frac{C\mu^{\frac{N-2}{2}}}{(1+\mu^2|y-{x}_{1}^{+}|^2)^\frac{N-2-\alpha}{2} }\left(\sum_{j=2}^{k}\frac{1}{\mu^\alpha|{x}_{j}^{+}-{x}_{1}^{+}|^\alpha}+\sum_{j=1}^{k}\frac{1}{\mu^\alpha|{x}_{1}^{+}-{x}_{j}^{-}|^\alpha}\right).
        \end{aligned}
    \end{equation}
Thus, we can obtain
\begin{align*}
    &\int_{\Omega_{1}^{+}}U_{{x}_{1}^{+}, \mu}^{2^*-2}\left(\sum_{j=2}^{k}U_{{x}_{j}^{+}, \mu}+\sum_{j=1}^{k}U_{{x}_{j}^{-},\mu}\right)^2\\
    \le& C\int_{\Omega_{1}^{+}} \frac{\mu^N}{(1+\mu^2
|y-{x}_{1}^{+}|^2)^{\frac{N+\sigma}{2}}}\left(\sum_{j=2}^{k}\frac{1}{\mu^\frac{N-\sigma}{2}|{x}_{j}^{+}-{x}_{1}^{+}|^\frac{N-\sigma}{2}}+\sum_{j=1}^{k}\frac{1}{\mu^\frac{N-\sigma}{2}|{x}_{1}^{+}-{x}_{j}^{-}|^\frac{N-\sigma}{2}}\right)^2\\
\le& \frac{C}{\mu^{N-\sigma}}\left(\sum_{j=2}^{k}\frac{1}{|{x}_{j}^{+}-{x}_{1}^{+}|^\frac{N-\sigma}{2}}+\sum_{j=1}^{k}\frac{1}{|{x}_{1}^{+}-{x}_{j}^{-}|^\frac{N-\sigma}{2}}\right)^2.
\end{align*}
Similarly, we have
\begin{align*}
&\int_{\Omega_{1}^{+}}\left(\sum_{j=2}^{k}U_{{x}_{j}^{+}, \mu}+\sum_{j=1}^{k}U_{{x}_{j}^{-},\mu}\right)^{2^*}
\le \frac{C}{\mu^{N-\sigma}}\left(\sum_{j=2}^{k}\frac{1}{|{x}_{j}^{+}-{x}_{1}^{+}|^\frac{(N-\sigma)(N-2)}{2N}}+\sum_{j=1}^{k}\frac{1}{|{x}_{1}^{+}-{x}_{j}^{-}|^\frac{(N-\sigma)(N-2)}{2N}}\right)^{2^*}.
\end{align*}
Therefore, we get
\begin{equation}\label{I4}
    \begin{aligned}
        I_4=&2k\left(\int_{\R^N}U_{0, 1}^{2^*}+2^* \left(\frac{B_0}{\mu^{N-2}}\left(\sum_{j=2}^{k}\frac{1}{|{x}_{j}^{+}-{x}_{1}^{+}|^{N-2}}+\sum_{j=1}^{k}\frac{1}{|{x}_{1}^{+}-{x}_{j}^{-}|^{N-2}}\right)\right)\right)
        \\&+k O\left(\frac{1}{\mu^N d_0^N}\right)+\frac{k}{\mu^{N-\sigma}}O\left(\sum_{j=2}^{k}\frac 1 {|{x}_{j}^{+}-{x}_{1}^{+}|^{N-\sigma}}+\sum_{j=1}^{k}\frac 1 {|{x}_{1}^{+}-{x}_{j}^{-}|^{N-\sigma}}\right)\\
&+\frac{k}{\mu^{N-\sigma}}O \left(\sum_{j=2}^{k}\frac{1}{|{x}_{j}^{+}-{x}_{1}^{+}|^\frac{N-\sigma}{2}}+\sum_{j=1}^{k}\frac{1}{|{x}_{1}^{+}-{x}_{j}^{-}|^\frac{N-\sigma}{2}}\right)^2
\\&+\frac{k}{\mu^{N-\sigma}}O\left(\sum_{j=2}^{k}\frac{1}{|{x}_{j}^{+}-{x}_{1}^{+}|^\frac{(N-\sigma)(N-2)}{2N}}+\sum_{j=1}^{k}\frac{1}{|{x}_{1}^{+}-{x}_{j}^{-}|^\frac{(N-\sigma)(N-2)}{2N}}\right)^{2^*}.
    \end{aligned}
\end{equation}

Combining \eqref{I}--\eqref{I2}  and \eqref{I4}, we have
\begin{equation*}
    \begin{aligned}
        I(W_{r,h,\mu})=&k(I_0+I_1+I_2+I_3)-\frac{1}{2^*}I_4\\
        =&k\left(\frac{2}{N} \int_{\R^N}U_{0, 1}^{2^*}+\frac{V(r)}{\mu^2} \int_{\R^N}U_{0, 1}^{2} -\frac{B_0}{\mu^{N-2}}\left(\sum_{j=2}^{k}\frac{1}{|{x}_{j}^{+}-{x}_{1}^{+}|^{N-2}}+\sum_{j=1}^{k}\frac{1}{|{x}_{1}^{+}-{x}_{j}^{-}|^{N-2}}\right)\right)\\
        &+ko\left(\frac{1}{\mu^2}\right)+k O\left(\frac{1}{\mu^N d_0^N}\right)+kO\left( \frac{1}{\mu^{N-2}}\left(\sum_{j=2}^{k}\frac 1 {|{x}_{j}^{+}-{x}_{1}^{+}|^{N-4}}+\sum_{j=1}^{k}\frac 1 {|{x}_{1}^{+}-{x}_{j}^{-}|^{N-4}}\right)  \right)\\
        &+\frac{k}{\mu^{N-\sigma}}O\left(\sum_{j=2}^{k}\frac 1 {|{x}_{j}^{+}-{x}_{1}^{+}|^{N-\sigma}}+\sum_{j=1}^{k}\frac 1 {|{x}_{1}^{+}-{x}_{j}^{-}|^{N-\sigma}}\right)\\
&+\frac{k}{\mu^{N-\sigma}}O \left(\sum_{j=2}^{k}\frac{1}{|{x}_{j}^{+}-{x}_{1}^{+}|^\frac{N-\sigma}{2}}+\sum_{j=1}^{k}\frac{1}{|{x}_{1}^{+}-{x}_{j}^{-}|^\frac{N-\sigma}{2}}\right)^2
\\&+\frac{k}{\mu^{N-\sigma}}O\left(\sum_{j=2}^{k}\frac{1}{|{x}_{j}^{+}-{x}_{1}^{+}|^\frac{(N-\sigma)(N-2)}{2N}}+\sum_{j=1}^{k}\frac{1}{|{x}_{1}^{+}-{x}_{j}^{-}|^\frac{(N-\sigma)(N-2)}{2N}}\right)^{2^*},
    \end{aligned}
\end{equation*}
    thus we have proved \eqref{I(W)}.
\end{proof}

\medskip
\begin{proposition}\label{A.4}
  Suppose that $V(|y|)$ satisfies $({\bf V})$ and $N\ge 5, (r,h,\mu)\in \mathscr{S}_k$, then we have the following expansion of the energy functional for $k \to +\infty$:
\begin{equation}\label{II(W)}
\begin{aligned}
    I(W_{r,h,\mu})=&k\left(A_1+\frac{A_2 V(r)}{\mu^2}-\frac{B_3 k^{N-2}}{\mu^{N-2}(\sqrt{1-h^2})^{N-2}}- \frac{B_4 k }{\mu^{N-2}h^{N-3}\sqrt{1-h^2}} \right)\\
    &+ko\left(\frac{1}{\mu^2}\right)+k O\left(\frac{k^{N-\sigma}}{\mu^{N-\sigma}}+\frac{k^{N-4} \ln k}{\mu^{N-2}}\right)+k O\left(\frac{1}{\mu^{N-2}h^{N-2}}\right),
\end{aligned}
\end{equation}
where $B_3=\displaystyle\frac{B_0B_1}{r^{N-2}}, B_4=\frac{B_0B_2}{r^{N-2}}$ are positive constants and  $ \sigma $ is a  constant small enough.
\end{proposition}
\begin{proof}
Using Lemma \ref{A.3} and Lemma \ref{B.1}, we can  obtain \eqref{II(W)}.
\end{proof}

\begin{remark}\label{A.5}
    We should first discuss how $h$ behaves as $k \to +\infty$.  Under the assumption of Lemma \ref{A.3}, if $h \neq o_k(1)$, then from Lemma \ref{B.1}, we can calculate the main term of the energy functional is $$k\displaystyle\left(A_1+\frac{A_2 V(r)}{\mu^2}-\frac{B_3k^{N-2}}{\mu^{N-2}(\sqrt{1-h^2})^{N-2}}\right).$$  Based on this, if $\sqrt{1-h^2} \neq o_k(1)$ neither, it is quite similar to the situation in \cite{WY}. In this case, we can construct bubble solutions of bilayer construction by doubling the bubble solutions in \cite{WY}. This kind of arguments are so standard that we omit it. Else if $\sqrt{1-h^2} = o_k(1)$, we are not able to find a critical point to the main term of the energy functional. That is why we may assume that $h = o_k(1)$, and additionally  supplement the assumption as $(r,h,\mu)\in \mathscr{S}_k$.
\end{remark}

\begin{proposition}\label{A.6}
    Suppose that $V(|y|)$ satisfies $({\bf V})$ and $N\ge 5, (r,h,\mu)\in \mathscr{S}_k$, then we have the following expansion for $k \to +\infty$:
\begin{equation}\label{Ih(W)}
\begin{aligned}
    \frac{\partial I(W_{r,h,\mu})}{\partial h}=&k\left(-\frac{(N-2)B_3 k^{N-2}h}{\mu^{N-2}(\sqrt{1-h^2})^{N}}+\frac{(N-3)B_4 k}{\mu^{N-2}h^{N-2}\sqrt{1-h^2}}-\frac{B_4 k}{\mu^{N-2}h^{N-4} (\sqrt{1-h^2})^3 }\right)
\\&+kO\left(  \frac{k^{N-3}\ln k}{\mu^{N-2}}\right)+kO\left(  \frac{1}{\mu^{N-2}h^{N-1}}+\frac{hk^{N-\sigma}}{\mu^{N-\sigma}}\right),
\end{aligned}
\end{equation}
and
\begin{equation}\label{Imu(W)}
\begin{aligned}
    \frac{\partial I(W_{r,h,\mu})}{\partial \mu}=&k\left(-\frac{2 A_2 V(r)}{\mu^3}+\frac{(N-2)B_3k^{N-2}}{\mu^{N-1}(\sqrt{1-h^2})^{N-2}}+ \frac{(N-2) B_4 k}{\mu^{N-1}h^{N-3}\sqrt{1-h^2}} \right)\\
    &+k O\left(\frac{k^{N-\sigma}}{\mu^{N+1-\sigma}}+\frac{k^{N-3}\ln k}{\mu^{N-1}}\right)+k O\left(\frac{1}{\mu^{N-1}h^{N-2}}\right),
\end{aligned}
\end{equation}
where $B_3, B_4$ are positive constants in Lemma \ref{A.4} and  $ \sigma $ is a  constant small enough.
\end{proposition}
\begin{proof}
Let  $\overline{\mathbb{Z}}_{\ell j}$ and $\underline{\mathbb{Z}}_{\ell j}$ be as following:
		\begin{align*}
			&\overline{\mathbb{Z}}_{1j}=\frac{\partial U_{x^{+}_{j}, \mu}}{\partial r},
			\qquad \qquad
			\overline{\mathbb{Z}}_{2j}=\frac{\partial U_{x^{+}_{j}, \mu}}{\partial h},
			\qquad \qquad
			\overline{\mathbb{Z}}_{3j}=\frac{\partial U_{x^{+}_{j}, \mu}}{\partial \mu}, \nonumber
			\\
			& \underline{\mathbb{Z}}_{1j}=\frac{\partial U_{x^{-}_{j}, \mu}}{\partial r},
			\qquad \qquad
			\underline{\mathbb{Z}}_{2j}=\frac{\partial U_{x^{-}_{j}, \mu}}{\partial h},
			\qquad \qquad
			\underline{\mathbb{Z}}_{3j}=\frac{\partial U_{x^{-}_{j}, \mu}}{\partial \mu},
		\end{align*}
		for $ j= 1, \cdots, k $.
Note that
   \begin{equation}\label{Idh}
       \begin{aligned}
     \frac{\partial I(W_{r,h,\mu})}{\partial h}=&2k\int_{\R^N} U_{{x}_{1}^{+}, \mu}^{2^*-1} \frac{\partial W_{r,h,\mu}}{\partial h}+\int_{\R^N} V(|y|) W_{r,h,\mu} \frac{\partial W_{r,h,\mu}}{\partial h}
     -\int_{\R^N} W_{r,h,\mu}^{2^*-1} \frac{\partial W_{r,h,\mu}}{\partial h}\\
     =&2k\int_{\R^N} U_{{x}_{1}^{+}, \mu}^{2^*-1}\left(\sum_{j=2}^k \overline{\mathbb Z}_{2j} +\sum_{j=1}^k \underline{\mathbb Z}_{2j}\right)+2k \int_{\R^N} V(|y|) U_{{x}_{1}^{+}, \mu}\overline{\mathbb Z}_{21}\\&+2k\int_{\R^N} V(|y|) U_{{x}_{1}^{+}, \mu}\left(\sum_{j=2}^k \overline{\mathbb Z}_{2j} +\sum_{j=1}^k \underline{\mathbb Z}_{2j}\right)-\int_{\R^N} W_{r,h,\mu}^{2^*-1} \frac{\partial W_{r,h,\mu}}{\partial h}.
\end{aligned}
   \end{equation}
And
\begin{equation}
    \begin{aligned}
&\int_{\R^N} W_{r,h,\mu}^{2^*-1} \frac{\partial W_{r,h,\mu}}{\partial h}\\
=&4k\int_{\R^N}U_{{x}_{1}^{+}, \mu}^{2^*-1}\left(\sum_{j=2}^k \overline{\mathbb Z}_{2j} +\sum_{j=1}^k \underline{\mathbb Z}_{2j}\right) +4k\int_{\R^N \setminus \Omega_{1}^{+}} U_{{x}_{1}^{+}, \mu}^{2^*-1} \frac{\partial W_{r,h,\mu}}{\partial h} \\
&+kO\left(\int_{\Omega_{1}^{+}}U_{{x}_{1}^{+}, \mu}^{2^*-2}\left(\sum_{j=2}^{k}U_{{x}_{j}^{+}, \mu}+\sum_{j=1}^{k}U_{{x}_{j}^{-}, \mu}\right) \left(\sum_{j=2}^k \overline{\mathbb Z}_{2j} +\sum_{j=1}^k \underline{\mathbb Z}_{2j}\right)\right)\\
&+kO\left(\int_{\Omega_{1}^{+}}\left(\sum_{j=2}^{k}U_{{x}_{j}^{+}, \mu}+\sum_{j=1}^{k}U_{{x}_{j}^{-}, \mu}\right)^{2^*-1}\frac{\partial W_{r,h,\mu}}{\partial h}\right)
    \end{aligned}
\end{equation}
Since for $j=1,2,\cdots,k,$
\begin{align*}
    |\overline{\mathbb Z}_{2j}|=&\left|\frac{rh\cos{\frac{2(j-1)\pi}{k}}}{\sqrt{1-h^2}}\frac{\partial U_{x^{+}_{j},\mu}}{\partial y_1}+\frac{rh\sin{\frac{2(j-1)\pi}{k}}}{\sqrt{1-h^2}}\frac{\partial U_{x^{+}_{j},\mu}}{\partial y_2}-r\frac{\partial U_{x^{+}_{j},\mu}}{\partial y_3}\right| \\
    =&\left|\frac{\partial }{\partial h}\frac{C_N \mu^{\frac{N-2}{2}}}{\left(1+\mu^2|y-{x}_{j}^{+}|^2\right)^{\frac{N-2}{2}}}\right|=\left|-(N-2)\frac{C_N \mu^{\frac{N}{2}}}{\left(1+\mu^2|y-{x}_{j}^{+}|^2\right)^{\frac{N}{2}}}\mu|y-{x}_{j}^{+}| \frac{\partial |y-{x}_{j}^{+}|}{\partial h}\right|\\
    =&\left|-(N-2)\frac{C_N \mu^{\frac{N}{2}}}{\left(1+\mu^2|y-{x}_{j}^{+}|^2\right)^{\frac{N}{2}}} \left\langle \mu(y-{x}_{j}^{+}),\left( \frac{rh\cos{\frac{2(j-1)\pi}{k}}}{\sqrt{1-h^2}} ,\frac{rh\sin{\frac{2(j-1)\pi}{k}}}{\sqrt{1-h^2}},-r,{\bf 0}\right) \right\rangle\right|\\
    \le&C\frac{ r\mu^{\frac{N}{2}}}{\left(1+\mu^2|y-{x}_{j}^{+}|^2\right)^{\frac{N-1}{2}}}.
\end{align*}
Meanwhile,
    \begin{align*}
    |\underline{\mathbb Z}_{2j}|=&\left|\frac{rh\cos{\frac{2(j-1)\pi}{k}}}{\sqrt{1-h^2}}\frac{\partial U_{x^{-}_{j},\mu}}{\partial y_1}+\frac{rh\sin{\frac{2(j-1)\pi}{k}}}{\sqrt{1-h^2}}\frac{\partial U_{x^{-}_{j},\mu}}{\partial y_2}+r\frac{\partial U_{x^{-}_{j},\mu}}{\partial y_3}\right| \\
    =&\left|-(N-2)\frac{C_N \mu^{\frac{N}{2}}}{\left(1+\mu^2|y-{x}_{j}^{-}|^2\right)^{\frac{N}{2}}}\left\langle \mu (y-{x}_{j}^{-}),\left( \frac{rh\cos{\frac{2(j-1)\pi}{k}}}{\sqrt{1-h^2}} ,\frac{rh\sin{\frac{2(j-1)\pi}{k}}}{\sqrt{1-h^2}},r,{\bf 0}\right) \right\rangle\right|\\
    \le&C\frac{ r\mu^{\frac{N}{2}}}{\left(1+\mu^2|y-{x}_{j}^{-}|^2\right)^{\frac{N-1}{2}}},
\end{align*}
 then by similar calculations in the proof of Lemma \ref{A.2} and Proposition \ref{A.4}, we can get by symmetry that for any small fixed $\delta>0$, it holds
 \begin{equation}\label{IV}
     \begin{aligned}
   & \int_{\R^N} V(|y|)U_{{x}_{1}^{+}, \mu}\overline{\mathbb Z}_{21} \\
    =&\frac{1}{\mu}\left(\int_{B_ {\mu\delta}(0)}\left(V(r)+V'(\xi(|\mu_k^{-1}y+{x}_{1}^{+}|))\left(|\mu_k^{-1}y+{x}_{1}^{+}|-r\right)\right)U_{0, 1} \overline{\mathbb Z}_{21} |_{x=0,\mu=1} \right)+O\left(\frac{1}{\mu^{N-2}}\right)\\
    =&O\left(\frac{1}{\mu^{N-2}}\right),
\end{aligned}
 \end{equation}
 the last equation holds because $\xi(|\mu_k^{-1}y+{x}_{1}^{+}|)$ is radial about ${x}_{1}^{+}.$ And,
\begin{equation}\label{IVI}
    \begin{aligned}
    &\int_{\R^N} V(|y|) U_{{x}_{1}^{+}, \mu}\left(\sum_{j=2}^k \overline{\mathbb Z}_{2j} +\sum_{j=1}^k \underline{\mathbb Z}_{2j}\right)\\
    =&O\left(\frac{1}{\mu^{N-2}}\left(\sum_{j=2}^{k}\frac{1}{|{x}_{j}^{+}-{x}_{1}^{+}|^{N-3}}+\sum_{j=1}^{k}\frac{1}{|{x}_{1}^{+}-{x}_{j}^{-}|^{N-3}}\right) \right)\\
    =&O\left( \frac{k^{N-3} \ln k}{\mu^{N-2}} \right).
\end{aligned}
\end{equation}
Meanwhile, for $j=2,\cdots,k,$ we can obtain
\begin{equation}
    \begin{aligned}
        &\int_{\R^N} U_{{x}_{1}^{+}, \mu}^{2^*-1} \overline{\mathbb Z}_{2j} \\ =&\int_{B_{\frac{\overline{{d}_j} }{2}}({x}_{1}^{+})}  U_{{x}_{1}^{+}, \mu}^{2^*-1} \frac{\partial U_{{x}_{j}^{+}, \mu}}{\partial h}+\int_{B_{\frac{\overline{{d}_j} }{2}}({x}_{j}^{+})}  U_{{x}_{1}^{+}, \mu}^{2^*-1} \frac{\partial U_{{x}_{j}^{+}, \mu}}{\partial h}+\int_{B^C_{ij}}  U_{{x}_{1}^{+}, \mu}^{2^*-1} \frac{\partial U_{{x}_{j}^{+}, \mu}}{\partial h}\\
 :=&I_1+I_2+I_3.
    \end{aligned}
\end{equation}
By Taylor expansion, we calculate for  $I_1$:
\begin{align}\label{hh1}
    I_1=&\int_{B_{\frac{\overline{{d}_j} }{2}}({x}_{1}^{+})}  U_{{x}_{1}^{+}, \mu}^{2^*-1} \frac{\partial U_{{x}_{j}^{+}, \mu}}{\partial h} \nonumber \\
    =&- \frac{r\mu(N-2)C_N^{2^*}}{\sqrt{1-h^2}}\int_{B_{\frac{\mu\overline{{d}_j} }{2}}(0)}\frac{1}{(1+z^{2})^{\frac{N+2} 2}}\left(  \frac{1}{\mu^{N}|{x}_{j}^{+}-{x}_{1}^{+}|^{N}}+O\left(\frac{1+z^2}{\mu^{N+2}|{x}_{j}^{+}-{x}_{1}^{+}|^{N+2}}\right) \right)\nonumber  \\
     &\times \left \langle \mu({x}_{1}^{+}-{x}_{j}^{+}),\left( h\cos{\frac{2(j-1)\pi}{k}},h\sin{\frac{2(j-1)\pi}{k}},-{\sqrt{1-h^2}},{\bf 0}\right) \right\rangle \text{d} z\nonumber  \\
   &+ \frac{r\mu N(N-2)C_N^{2^*}}{\sqrt{1-h^2}}\int_{B_{\frac{\mu\overline{{d}_j} }{2}}(0)}\frac{1}{(1+z^{2})^{\frac{N+2} 2}}  \left(\frac{ \langle z ,\mu ({x}_{1}^{+}-{x}_{j}^{+}) \rangle}{\mu^{N+2}|{x}_{j}^{+}-{x}_{1}^{+}|^{N+2}}+O\left(\frac{(1+z^{2}+2  \mu\langle z , {x}_{1}^{+}-{x}_{j}^{+} \rangle)^{2}}{\mu^{N+4}|{x}_{j}^{+}-{x}_{1}^{+}|^{N+4}}\right) \right)\nonumber \\
   &\times \left \langle z,\left( h\cos{\frac{2(j-1)\pi}{k}},h\sin{\frac{2(j-1)\pi}{k}},-{\sqrt{1-h^2}},{\bf 0}\right) \right\rangle \text{d} z\nonumber  \\
=&\frac{B_0 r^2 h (N-2)}{\mu^{N-2}|{x}_{j}^{+}-{x}_{1}^{+}|^{N}}\left(1-\cos{\frac{2(j-1)\pi}{k}}\right)+O\left(\frac{h\left(1-\cos{\frac{2(j-1)\pi}{k}}\right) }{\mu^{N-\sigma}|{x}_{j}^{+}-{x}_{1}^{+}|^{N+2-\sigma}} \right).
\end{align}
Next, we calculate $I_2,$ we have
\begin{equation}\label{hh2}
    \begin{aligned}
 I_2=&\int_{B_{\frac{\overline{{d}_j} }{2}}({x}_{j}^{+})}  U_{{x}_{1}^{+}, \mu}^{2^*-1} \frac{\partial U_{{x}_{j}^{+}, \mu}}{\partial h}\\
 =&- \frac{r\mu(N-2)C_N^{2^*}}{\sqrt{1-h^2}}\int_{B_{\frac{\mu\overline{{d}_j} }{2}}(0)}  \frac{1}{(1+z^{2})^{\frac{N} 2}}  \frac{1}{(1+|z+\mu{x}_{j}^{+}-\mu{x}_{1}^{+}|^{2})^{\frac{N+2}{2}}} \\
    &\times \left \langle z,\left( h\cos{\frac{2(j-1)\pi}{k}},h\sin{\frac{2(j-1)\pi}{k}},-{\sqrt{1-h^2}},{\bf 0}\right) \right\rangle \text{d} z 
    =O\left(\frac{h\left(1-\cos{\frac{2(j-1)\pi}{k}}\right)}{\mu^{N}|{x}_{j}^{+}-{x}_{1}^{+}|^{N+2}}\right).
\end{aligned}
\end{equation}
As for $I_3$, by a standard estimate, we have
\begin{equation*}
    \begin{aligned}
 |I_3|\le \int_{\overline{B^C_{1j}}} \left| U_{{x}_{1}^{+}, \mu}^{2^*-1} \frac{\partial U_{{x}_{j}^{+}, \mu}}{\partial h}\right|=O\left(\frac{1}{\mu^{N}|{x}_{j}^{+}-{x}_{1}^{+}|^{N+1}\sqrt{1-h^2}}\right).
 \end{aligned}
\end{equation*}
However, compared to the main term of $I_1$, this estimate is too rough. So we must obtain a more accurate estimate. By a symmetry trick, we can get:
\begin{equation}\label{hh3}
    \begin{aligned}
 I_3=\int_{\overline{B^C_{1j}}}  U_{{x}_{1}^{+}, \mu}^{2^*-1} \frac{\partial U_{{x}_{j}^{+}, \mu}}{\partial h}=O\left(\frac{h\sin {\frac{2(j-1)\pi}{k}}}{\mu^{N}|{x}_{j}^{+}-{x}_{1}^{+}|^{N+1}\sqrt{1-h^2}}\right).
 \end{aligned}
\end{equation}

From \eqref{hh1}--\eqref{hh3}, we have
\begin{equation*}\label{hh+}
    \begin{aligned}
 &\sum_{j=2}^{k}\int_{\R^N} U_{{x}_{1}^{+}, \mu}^{2^*-1} \overline{\mathbb Z}_{2j}\\
 =&\sum_{j=2}^{k} \frac{B_0 r^2 h (N-2)}{\mu^{N-2}|{x}_{j}^{+}-{x}_{1}^{+}|^{N}}\left(1-\cos{\frac{2(j-1)\pi}{k}}\right)\\
 &+O\left(\sum_{j=2}^{k}\frac{h\left(1-\cos{\frac{2(j-1)\pi}{k}}\right)\ln(\mu |{x}_{j}^{+}-{x}_{1}^{+}|) }{\mu^{N}|{x}_{j}^{+}-{x}_{1}^{+}|^{N+2}} \right)+ O\left(\sum_{j=2}^{k}\frac{h\sin {\frac{2(j-1)\pi}{k}}}{\mu^{N}|{x}_{j}^{+}-{x}_{1}^{+}|^{N+1}\sqrt{1-h^2}}\right)\\
  =&\sum_{j=2}^{k} \frac{B_0 r^2 h (N-2)}{\mu^{N-2}|{x}_{j}^{+}-{x}_{1}^{+}|^{N}}\left(1-\cos{\frac{2(j-1)\pi}{k}}\right)+O\left(\frac{hk^{N-\sigma}}{\mu^{N-\sigma}}\right).
 \end{aligned}
\end{equation*}
By direct calculating, we also have
\begin{equation*}\label{hh-}
    \begin{aligned}
 \sum_{j=1}^{k}\int_{\R^N} U_{{x}_{1}^{+}, \mu}^{2^*-1} \underline{\mathbb Z}_{2j}=\sum_{j=1}^{k} \frac{B_0 r^2 h (N-2)}{\mu^{N-2}|{x}_{1}^{+}-{x}_{j}^{-}|^{N}}\left(1-\cos{\frac{2(j-1)\pi}{k}}\right)-\sum_{j=1}^{k}\frac{2 B_0 r^2 h (N-2)}{\mu^{N-2}|{x}_{1}^{+}-{x}_{j}^{-}|^{N}}+O\left(\frac{k}{\mu^{N-\sigma}h^{N-\sigma}}\right).
 \end{aligned}
\end{equation*}
Therefore, using Lemma \ref{B.1} and Lemma \ref{B.2}, we have
\begin{equation}\label{hhh}
    \begin{aligned}
 &\int_{\R^N} U_{{x}_{1}^{+}, \mu}^{2^*-1} \left(\sum_{j=2}^k \overline{\mathbb Z}_{2j} +\sum_{j=1}^k \underline{\mathbb Z}_{2j}\right)\\
  =&\left(\sum_{j=2}^{k} \frac{B_0 r^2 h (N-2)}{\mu^{N-2}|{x}_{j}^{+}-{x}_{1}^{+}|^{N}}\left(1-\cos{\frac{2(j-1)\pi}{k}}\right)+\sum_{j=1}^{k} \frac{B_0 r^2 h (N-2)}{\mu^{N-2}|{x}_{1}^{+}-{x}_{j}^{-}|^{N}}\left(1-\cos{\frac{2(j-1)\pi}{k}}\right)\right)\\
  &-\sum_{j=1}^{k}\frac{2 B_0 r^2 h (N-2)}{\mu^{N-2}|{x}_{1}^{+}-{x}_{j}^{-}|^{N}}+O\left(\frac{hk^{N-\sigma}}{\mu^{N-\sigma}}\right)+O\left(\frac{k}{\mu^{N-\sigma}h^{N-\sigma}}\right).
 \end{aligned}
\end{equation}

On the other hand, we have
\begin{equation}\label{WWh}
    \begin{aligned}
&\int_{\R^N} W_{r,h,\mu}^{2^*-1} \frac{\partial W_{r,h,\mu}}{\partial h}\\
=&4k\int_{\R^N}U_{{x}_{1}^{+}, \mu}^{2^*-1}\left(\sum_{j=2}^k \overline{\mathbb Z}_{2j} +\sum_{j=1}^k \underline{\mathbb Z}_{2j}\right) \\
&+kO\left(\sum_{j=2}^{k} \frac{h\left(1-\cos{\frac{2(j-1)\pi}{k}}\right)}{\mu^{N-2}|{x}_{j}^{+}-{x}_{1}^{+}|^{N}}+\sum_{j=1}^{k} \frac{h\left(1-\cos{\frac{2(j-1)\pi}{k}}\right)}{\mu^{N-2}|{x}_{1}^{+}-{x}_{j}^{-}|^{N}}\right)\\ &\times\left(\frac{1}{\mu^{\frac{N-2}{2}}|{x}_{2}^{+}-{x}_{1}^{+}|^{\frac{N-2}{2}}}+\frac{1}{\mu^{\frac{N-2}{2}}|{x}_{1}^{+}-{x}_{1}^{-}|^{\frac{N-2}{2}}}\right)\\
=&4k\int_{\R^N}U_{{x}_{1}^{+}, \mu}^{2^*-1}\left(\sum_{j=2}^k \overline{\mathbb Z}_{2j} +\sum_{j=1}^k \underline{\mathbb Z}_{2j}\right) +kO\left(\frac{hk^{N-2}}{\mu^{N-1}} \right).
    \end{aligned}
\end{equation}
Then from \eqref{Idh}, \eqref{IV}, \eqref{IVI}, \eqref{hhh}, \eqref{WWh} and Lemma \ref{B.2}, we finally get
\begin{equation}
    \begin{aligned}
 &\frac{\partial I(W_{r,h,\mu})}{\partial h}\\
     =&-2k\int_{\R^N} U_{{x}_{1}^{+}, \mu}^{2^*-1}\left(\sum_{j=2}^k \overline{\mathbb Z}_{2j} +\sum_{j=1}^k \underline{\mathbb Z}_{2j}\right)+2k \int_{\R^N} V(|y|) U_{{x}_{1}^{+}, \mu}\overline{\mathbb Z}_{21}\\
     &+2k\int_{\R^N} V(|y|) U_{{x}_{1}^{+}, \mu}\left(\sum_{j=2}^k \overline{\mathbb Z}_{2j} +\sum_{j=1}^k \underline{\mathbb Z}_{2j}\right)+kO\left(\frac{hk^{N-\sigma}}{\mu^{N-\sigma}}\right)+kO\left(\frac{k}{\mu^{N-\sigma}h^{N-\sigma}}\right)+kO\left(\frac{hk^{N-2}}{\mu^{N-1}} \right)\\
=&k\left(-\frac{(N-2)B_3 k^{N-2}h}{\mu^{N-2}(\sqrt{1-h^2})^{N}}+\frac{(N-3)B_4 k}{\mu^{N-2}h^{N-2}\sqrt{1-h^2}}-\frac{B_4 k}{\mu^{N-2}h^{N-4} (\sqrt{1-h^2})^3 }\right)
\\&+kO\left(  \frac{k^{N-3}\ln k}{\mu^{N-2}}\right)+kO\left(  \frac{1}{\mu^{N-2}h^{N-1}}\right)+kO\left(\frac{hk^{N-\sigma}}{\mu^{N-\sigma}}\right).
    \end{aligned}
\end{equation}

Similarly, we can also prove \eqref{2*-1+-}. Noting that for $j=1,2,\cdots,k,$
\begin{align*}
    |\overline{\mathbb Z}_{3j}|=&\left|\frac{N-2}{2}\frac{C_N \mu^{\frac{N-4}{2}}}{\left(1+\mu^2|y-{x}_{j}^{+}|^2\right)^{\frac{N-2}{2}}}-\frac{N-2}{2}\frac{2C_N \mu^{\frac{N}{2}}|y-{x}_{j}^{+}|^2}{\left(1+\mu^2|y-{x}_{j}^{+}|^2\right)^{\frac{N}{2}}}\right|\\
    =&\left|\frac{N-2}{2}\frac{C_N \mu^{\frac{N-4}{2}}}{\left(1+\mu^2|y-{x}_{j}^{+}|^2\right)^{\frac{N-2}{2}}}\frac{1-\mu^2|y-{x}_{j}^{+}|^2}{1+\mu^2|y-{x}_{j}^{+}|^2}\right|
    \le\frac{U_{{x}_{j}^{+}, \mu}}{\mu},
\end{align*}
and
\begin{align*}
    |\underline{\mathbb Z}_{3j}|
    =\left|\frac{N-2}{2}\frac{C_N \mu^{\frac{N-4}{2}}}{\left(1+\mu^2|y-{x}_{j}^{-}|^2\right)^{\frac{N-2}{2}}}\frac{1-\mu^2|y-{x}_{j}^{+}|^2}{1+\mu^2|y-{x}_{j}^{-}|^2}\right|
    \le\frac{U_{{x}_{j}^{-}, \mu}}{\mu}.
\end{align*}
By some calculations based on Proposition \ref{pro2.2}, we can obtain \eqref{Imu(W)}.

\end{proof}
\medskip

	\section{Finite dimensional reduction}  \label{sec2}
In this section, we perform a finite-dimensional reduction. Set $\mathbb{E}_k$  be given by
		\begin{align}\label{SpaceE}
			\mathbb{E}_k= \Big\{\phi:  \phi\in H_s, \quad \int_{\R^N} U_{x^{+}_{j}, \mu}^{2^*-2}  \overline{\mathbb{Z}}_{\ell j}\phi = \int_{\R^N} U_{x^{+}_{j}, \mu}^{2^*-2}  \underline{\mathbb{Z}}_{\ell j}\phi = 0,  \quad j= 1, \cdots, k, \quad \ell= 1,2, 3  \Big\}.
            \end{align}

	We consider the following linearized problem:
	\begin{align}\label{lin}
		\begin{cases}
			{\bf L_k}\phi_k={\bf N}(\phi_k)+{\bf l_k}+\sum\limits_{i=1}^k\sum\limits_{\ell=1}^3 \Big({c}_\ell \overline{\mathbb{Z}}_{\ell i}+{c}_\ell \underline{\mathbb{Z}}_{\ell i}\Big)\, \; \text{in}\;
			\R^N,\\
			\phi_k \in  \mathbb{E}_k,
		\end{cases}
	\end{align}
	for some constants $c_{\ell} $, where
 \begin{equation*}
		\langle{{\bf L_k}}\phi_k,\omega\rangle=\int_{\R^N}\left(-\Delta {\phi_k}+V(|y|)\phi_k-(2^*-1)W_{r,h,\mu}^{2^*-2}\phi_k\right)\omega,
	\end{equation*}
	   \begin{equation*}
		\langle{{\bf N}}(\phi_k),\omega\rangle=\int_{\R^N}\left(\bigl(W_{r,h,\mu}+\phi_k \bigr)^{2^*-1}-W_{r,h,\mu}^{2^*-1}-(2^*-1)W_{r,h,\mu}^{2^*-2}\phi_k\right)\omega,
	\end{equation*}
	and
 \begin{equation*}
		\langle{\bf l}_k,\omega\rangle=\int_{\R^N}\left(W_{r,h,\mu}^{2^*-1}-\sum_{j=1}^k \left(U_{x^{+}_j,\mu}^{2^*-1}+U_{x^{-}_j,\mu}^{2^*-1} \right)-V(|y|)W_{r,h,\mu}\right)\omega.
	\end{equation*}	
 We also have
  \begin{equation*}
      \langle\overline{\mathbb{Z}}_{\ell j},\omega\rangle=\int_{\R^N}\nabla \overline{\mathbb{Z}}_{\ell j}\nabla \omega=(2^*-1)\int_{\R^N}U_{x^{+}_{j}, \mu}^{2^*-2}  \overline{\mathbb{Z}}_{\ell j}\omega,  \quad\langle\underline{\mathbb{Z}}_{\ell j},\omega\rangle=(2^*-1)\int_{\R^N}U_{x^{-}_{j}, \mu}^{2^*-2}  \underline{\mathbb{Z}}_{\ell j}\omega.
  \end{equation*}

	\medskip
	
	\begin{lemma}\label{lem2.1}
		Suppose that $V(|y|)$ satisfies $({\bf V})$, $N\ge 5$ and $(r_k,h_k,\mu_k)\in {\mathscr S_k}$, there exists $\rho>0$ and $k_0>0$, such that for any $k\ge k_0$,
  \begin{equation}
      ||{\bf L_k}\phi||\ge \rho||\phi||,\qquad \forall \phi \in \mathbb{E}_k.
  \end{equation}
	\end{lemma}
 \begin{proof}
We prove by contradiction.  Suppose that  when  $k\rightarrow +\infty, (r_k, h_k, \mu_k)\in  \mathscr S_k $ and $\phi_k \in  \mathbb{E}_k$ satisfying
\begin{align}\label{o(1)}
\|{\bf{L_k}} \phi_k\| =  o_k(1) \| \phi_k\|.
\end{align}
Without loss of generality, we may assume that
\begin{equation*}
    ||\phi_k||^2=2k,\qquad\|{\bf{L_k}} \phi_k\|^2 =  o_k(k).
\end{equation*}
Then,
\begin{equation}\label{bound}
    \int_{\Omega_1^+} |D\phi_k|^2=1,
\end{equation}
and
\begin{equation}\label{o11}
    \int_{\Omega_1^+} \left(-\Delta {\phi_k}+V(|y|)\phi_k-(2^*-1)W_{r,h,\mu}^{2^*-2}\phi_k\right)\omega=o_k(1)\quad\forall \omega\in \mathbb{E}_k.
\end{equation}

We set $\widetilde{\phi_k}(y)=\mu_k^{-\frac{N-2}{2}}\phi_k(\mu_k^{-1}y+{x}_{1}^{+}),$ then from \eqref{bound} we know that $\widetilde{\phi_k}$ is bounded in $D^{1,2}(\R^N).$ In fact, for any $R>0$, we can choose $k$ large enough such that $B_{\frac{R}{k}}(x^{+}_1) \subset \Omega_1^+.$ As a result, we have
\begin{equation*}
    \begin{aligned}
  \int_{B_{\frac{\mu_k R}{k}}(0)}| D\widetilde{\phi_k} |^2  \le  \int_{B_{\frac{R}{k}}(x^{+}_1)}  | D \phi_k |^2 \le 1.
    \end{aligned}
\end{equation*}
So we can deduce that there is a $\phi \in D^{1,2}(\R^N)$ that
 \begin{equation*}
    \widetilde{\phi_k}\rightarrow   \phi \quad  \text{ weakly in } D^{1,2}(\R^N),
 \end{equation*}
 and  \begin{equation*}
      \widetilde{\phi_k}\rightarrow  \phi \quad  \text{ strongly in } L^2_{loc}(\mathbb{R}^N).
 \end{equation*}
Since ${\phi_k}$ is even in $y_j, j = 2, 4,  \cdots, N$, thus $ \phi$ is even in  $y_j$.

From the orthogonal conditions for functions of $\mathbb{E}_k,$ we have
\begin{align*}
\int_{{\mathbb{R}}^N }  U_{x^{+}_1,\mu}^{{2^*-2}}  \frac{\partial U_ {x^{+}_1,\mu} }{\partial r } \phi_k  = 0
\end{align*}
and the identity
\begin{align*}
\frac{\partial  U_ {x^{+}_1,\mu} }  {  \partial r } =   \sqrt{1-h^2}  \frac{ \partial U_ {x^{+}_1,\mu} }{\partial {y_1} } + h \frac{ \partial U_ {x^{+}_1,\mu} }{\partial {y_3} },
\end{align*}
hence
\begin{align} \label{y1}
  \sqrt{1-h^2}\int_{{\mathbb{R}}^N }  U_ {x^{+}_1,\mu}^{{2^*-2}} \frac{ \partial U_{x^{+}_1,\mu} }{\partial {y_1} } \phi_k +  h  \int_{{\mathbb{R}}^N }  U_ {x^{+}_1,\mu}^{{2^*-2}} \frac{ \partial U_{x^{+}_1,\mu} }{\partial {y_3} } \phi_k=0.
  \end{align}
Similarly, combining
\begin{align*}
\int_{{\mathbb{R}}^N }  U_{x^{+}_1,\mu}^{{2^*-2}}  \frac{\partial U_ {x^{+}_1,\mu} }{\partial h } \phi_k = 0,
\end{align*}
and
\begin{align*}
\frac{\partial U_ {x^{+}_1,\mu} }{\partial h } = -\frac{ h r}{ \sqrt{1-h^2} }  \frac{ \partial U_ {x^{+}_1,\mu} }{\partial {y_1} }+ r \frac{ \partial U_ {x^{+}_1,\mu} }{\partial {y_3} },
\end{align*}
we  get
\begin{align} \label{y3}
  \frac{ h }{ \sqrt{1-h^2} }    \int_{{\mathbb{R}}^N }  U_ {x^{+}_1,\mu}^{{2^*-2}}       \frac{ \partial U_{x^{+}_1,\mu} }{\partial {y_1} } \phi_k  - \int_{{\mathbb{R}}^N }  U_ {x^{+}_1,\mu}^{{2^*-2}}       \frac{ \partial U_{x^{+}_1,\mu} }{\partial {y_3} }   \phi_k=0.
\end{align}
It follows from \eqref{y1}, \eqref{y3},  we obtain
\begin{align*}
\int_{{\mathbb{R}}^N }  U_ {x^{+}_1,\mu}^{{2^*-2}}  \frac{\partial U_ {x^{+}_1,\mu}  }{\partial y_1 }  \phi_k  = \int_{{\mathbb{R}}^N }  U_ {x^{+}_1,\mu}^{{2^*-2}}  \frac{\partial U_ {x^{+}_1,\mu}  }{\partial y_3 }  \phi_k =0.
\end{align*}
Thus,
\begin{align*}
\int_{{\mathbb{R}}^N }  U_ {0,1}^{{2^*-2}}  \frac{\partial U_ {0,1}  }{\partial y_1 }  \widetilde{\phi_k}  = \int_{{\mathbb{R}}^N }  U_ {0,1}^{{2^*-2}}  \frac{\partial U_ {0,1}  }{\partial y_3 }  \widetilde{\phi_k} =0.
\end{align*}
Taking $k \rightarrow   + \infty$, we obtain
 \begin{align}\label{condition3}
 \int_{{\mathbb{R}}^N } U_ {0,1}^{{2^*-2}} \frac{ \partial U_ {0,1} }{\partial {y_1} }     \phi = \int_{{\mathbb{R}}^N }  U_ {0,1}^{{2^*-2}}  \frac{ \partial U_ {0,1} }{\partial {y_3} } \phi=0.
 \end{align}

We claim that
\begin{equation}\label{claim}
    -\Delta\phi-(2^*-1)U_{0,1}^{2^*-2}\phi=0 \quad \text{in} \ \R^N.
\end{equation}
To prove \eqref{claim}, we define the space
 \begin{align*}
  {E}^+ = \Bigg\{\phi : \phi \in \overline{H^1(\R^N) \cap D^{1,2}(\R^N)}, \int_{{\mathbb{R}}^N}   U_{0,1}^{2^*-2}       \frac{ \partial U_{0,1} }{\partial {y_1} }  \phi= \int_{{\mathbb{R}}^N }  U_{0,1}^{2^*-2}  \frac{ \partial U_{0,1} }{\partial {y_3} } \phi =0\Bigg\}.
  \end{align*}
We first prove that $\phi $ is a solution of
  \begin{align}\label{y13}
- \Delta \phi - (2^*-1)U_{0,1}^{2^*-2} \phi  =0 \qquad \text{in}~   {E}^+.
\end{align}
In fact, for any $R>0, \omega\in C_0^{\infty}(B_{R}(0)) \cap {E}^+,$ we denote
\begin{equation*}
    \omega_k(y):=\mu_k^{\frac{N-2}{2}}\omega(\mu_k y-{x}_{1}^{+})\in C_0^{\infty}(B_{\mu_k^{-1}R}({x}_{1}^{+})).
\end{equation*}
Then we can deduce from \eqref{o(1)} that
\begin{equation*}
    \int_{\Omega_1^+ } \left(D{\phi_k}D\omega_k+V(|y|)\phi_k\omega_k-(2^*-1)W_{r,h,\mu}^{2^*-2}\phi_k\omega_k\right)=o_k\left(\frac{1}{\sqrt{k}}\right)||\omega_k||.
\end{equation*}
That is,
\begin{equation*}
    \int_{\R^N} \left(D\widetilde{\phi_k}D\omega+\mu_k^{-2}V\widetilde{\phi_k}\omega-(2^*-1)\left(\sum_{j=1}^k \left(U_{x^{+}_j-{x}_{1}^{+},1}+U_{x^{-}_j-{x}_{1}^{+},1}\right)\right)^{2^*-2}\widetilde\phi_k\omega\right)=o_k\left(\frac{1}{\sqrt{k}}\right)||\omega||.
\end{equation*}
Taking $k \to +\infty$, by Lebesgue's convergence theorem and Lemma \ref{B.3}, we have
\begin{equation}
    \int_{\mathbb{R}^N } \left(D\phi D\omega-(2^*-1)U_{0,1}^{2^*-2}\phi\omega\right)=0 \qquad\forall \omega\in C_0^{\infty}(B_{ R}(0)) \cap {E}^+, R>0.
\end{equation}
From the density of $\displaystyle\cup_{R>0}(C_0^{\infty}(B_{R}(0))\cap {E}^+)$ on $ {E}^+$, we have proved \eqref{y13}. Combining \eqref{y13} and the fact that \eqref{claim} holds for $\displaystyle\phi = \frac{\partial U}{\partial y_1} $ and $ \displaystyle\phi = \frac{\partial U}{\partial y_3}$,  we get
\begin{equation*}
   \int_{\mathbb{R}^N }  D\phi D \omega- (2^*-1)U_{0,1}^{2^*-2} \phi  \omega = 0,   \quad \forall \omega \in  \overline{H^1(\R^N) \cap D^{1,2}(\R^N)}.
\end{equation*}
Thus, we have $\phi=0$. Therefore, for any $R>0,$
\begin{equation}\label{ball}
   \int_{B_{\frac{R}{\mu_k}}(x_1^{+})} W_{r_k,h_k,\mu_k}^{2^*-2}{\phi_k}^2 \le\int_{B_{R}(0)} |\widetilde {\phi_k}|^2=o_k(1).
\end{equation}
Meanwhile, it is obvious that
\begin{equation*}
     W_{r_k,h_k,\mu_k}^{2^*-2}=o_k(1)\qquad \text{in} ~ \Omega_1^+ \setminus B_{\frac{R}{\mu_k}}(0).
\end{equation*}
Thus
\begin{equation*}
\begin{aligned}
        &(2^*-1)\int_{\Omega_1^+ \setminus B_{\frac{R}{\mu_k}}(x_1^{+})} W_{r_k,h_k,\mu_k}^{2^*-2}{\phi_k}^2\\
        \le&o(1) \int_{\Omega_1^+ } {\phi_k}^{2}+ \left(\int_{\Omega_1^+ \setminus B_{\frac{R}{\mu_k}}(x_1^{+})\setminus B_{r+\frac{R}{\mu_k}}(0)}  W_{r_k,h_k,\mu_k}^{2^*} \right)^{\frac{2}{N}} \left(\int_{\Omega_1^+ } {\phi_k}^{2^*} \right)^{\frac{2}{2^*}}\\
        \le& o(1) \int_{\Omega_1^+ } V(|y|){\phi_k}^{2}+\left(\frac{Ck^{2^*}}{\mu_k^N}\right)^{\frac{2}{N}}k^\frac{2}{N}\int_{\Omega_1^+ } {|D\phi_k|}^{2}\\=&o_k(1)\int_{\Omega_1^+ } \left({|D\phi_k|}^{2}+V(|y|){\phi_k}^{2}\right).
\end{aligned}
\end{equation*}
As a result, we get
\begin{equation}\label{ball-}
    \begin{aligned}
        (2^*-1) \int_{\Omega_1^+} W_{r_k,h_k,\mu_k}^{2^*-2}{\phi_k}^2 \le& C\left( \int_{B_{\frac{R}{\mu_k}}(x_1^{+})} W_{r_k,h_k,\mu_k}^{2^*-2}{\phi_k}^2+\int_{\Omega_1^+ \setminus B_{\frac{R}{\mu_k}}(x_1^{+})} W_{r_k,h_k,\mu_k}^{2^*-2}{\phi_k}^2\right)\\
        =& o_k(1)\int_{\Omega_1^+ } \left({|D\phi_k|}^{2}+V(|y|){\phi_k}^{2}\right)+o_k(1)\\=& o_k(1)\int_{\Omega_1^+ } \left({|D\phi_k|}^{2}+V(|y|){\phi_k}^{2}\right).
    \end{aligned}
\end{equation}
Combining \eqref{o11}, \eqref{ball} and \eqref{ball-}, we get
\begin{align*}
    o_k(1)=&\int_{\Omega_1^+ }\left({|D\phi_k|}^{2}+V(|y|){\phi_k}^{2}-(2^*-1)W_{r_k,h_k,\mu_k}^{2^*-2}{\phi_k}^2\right)\\
    =&(1+o_k(1))\int_{\Omega_1^+ } \left({|D\phi_k|}^{2}+V(|y|){\phi_k}^{2}\right)=o_k(1).
\end{align*}
So
\begin{equation}
    \int_{\Omega_1^+ } {|D\phi_k|}^{2}\le \int_{\Omega_1^+ } \left({|D\phi_k|}^{2}+V(|y|){\phi_k}^{2}\right)=o_k(1),
\end{equation}
which is a contradiction to \eqref{bound}. The proof of  Lemma \ref{lem2.1} is finished.
\end{proof}

\medskip

Using Fredholm alternative, we can	easily prove the following proposition by standard arguments.
	\begin{proposition}\label{pro2.2}
		There exist  $k_0>0  $ and a constant  $C>0 $ such
		that $\bf L_k$ is an isomorphism in $\mathbb{E}_k $.
	\end{proposition}
	\medskip

	Next, we will use the Contraction Mapping Principle to show that problem \eqref{lin} has a unique solution if  $\|\phi_k\| $ is small enough. For this purpose, we will first give the estimates  of $ {\bf l}_k $  and  ${{\bf N}}(\phi_k) $.

	\medskip

	\begin{lemma} \label{lem2.4}
		Suppose that $V(|y|)$ satisfies $({\bf V})$, $N\ge 5$ and $(r_k,h_k,\mu_k)\in {\mathscr S_k}$, then
		\begin{align}  \label{lk}
			\|\, {\bf l}_k \,\|  \le  Ck\max\left\{ \frac{1}{k^{\frac{1}{2^*}}}  \left( \frac{k}{\mu} \right)^{\frac{N+2}{2}-\tau} , \frac{1}{\mu^{\min\{{\frac{N-2}{2}},{2-\tau}\}}}\right\}.
		\end{align}
		where  $\tau$ is a small constant.
	\end{lemma}

	\begin{proof}
		We can rewrite $\langle{\bf l}_k,\omega\rangle$ as
		\begin{equation}\label{lk'}
		    \langle{\bf l}_k,\omega\rangle=\int_{\R^N}\left(W_{r,h,\mu}^{2^*-1}-\sum_{j=1}^k \left(U_{x^{+}_j,\mu}^{2^*-1}+U_{x^{-}_j,\mu}^{2^*-1} \right)\right)\omega-\int_{\R^N}V(|y|)W_{r,h,\mu}\omega:=I_1-I_2.
		\end{equation}
		For $I_2$, we have
  \begin{equation}
      \begin{aligned}
    |I_2|=&2k\left|\int_{\R^N}V(|y|)U_{x^{+}_1,\mu}\omega\right|=2k\left(\left|\int_{B_1(x^{+}_1)}V(|y|)U_{x^{+}_1,\mu}\omega\right|+\left|\int_{\R^N\setminus B_1(x^{+}_1)}V(|y|)U_{x^{+}_1,\mu}\omega\right|\right)\\
\le&Ck\left(\int_{B_1(x^{+}_1)}V(|y|)U_{x^{+}_1,\mu}|\omega|+\frac{||\omega||}{\mu^{\frac{N-2}{2}}}   \right)\\
          =&kO\left(\frac{1}{\mu^{\min\{{\frac{N-2}{2}},{2-\tau}
        \}}}\right)||\omega||.
      \end{aligned}
  \end{equation}
   The last inequality holds for the reason that
   \begin{align*}
    \left|\int_{B_1(x^{+}_1)}V(|y|)U_{x^{+}_1,\mu}\omega\right|\le& C \frac{1}{\mu^{N-\frac{(N-2)N}{N+2}}}\left(\int_{B_\mu(0)}U_{0,1}^{\frac{2^*}{2^*-1}}\right)^{\frac{2^*-1}{2^*}}||\omega||\\
    =&\left\{ \begin{matrix}
        O\left(\displaystyle\frac{1}{\mu^{\frac{3}{2}}}\right)||\omega||, \quad N=5,\\
        \\
        O\left(\displaystyle\frac{(\ln \mu)^{\frac{2}{3}}}{\mu^{2}}\right)||\omega||, \quad N=6,\\
        \\
        O\left(\displaystyle\frac{1}{\mu^{2}} \right)||\omega||, \quad N\ge7, \end{matrix}\right.\\
        =&O\left(\frac{1}{\mu^{\min\{{\frac{N-2}{2}},{2-\tau}
        \}}}\right).
   \end{align*}

Next, we calculate $I_1$. Noting that
\begin{equation}
\begin{aligned}
|I_1|=&2k\left|\int_{\Omega^+_1}\left(W_{r,h,\mu}^{2^*-1}-\sum_{j=1}^k \left(U_{x^{+}_j,\mu}^{2^*-1}+U_{x^{-}_j,\mu}^{2^*-1} \right)\right)\omega\right|\\
 \le& Ck \int_{\Omega^+_1}U_{x^{+}_1,\mu}^{2^*-2}\left(\sum_{j=2}^k U_{x^{+}_j,\mu}+\sum_{j=1}^kU_{x^{-}_j,\mu}\right) \left|\omega\right|+Ck\int_{\Omega^+_1}\left(\sum_{j=2}^k U_{x^{+}_j,\mu}+\sum_{j=1}^kU_{x^{-}_j,\mu}\right)^{2^*-1}\left|\omega\right|\\
 \le& Ck(I_{11}+I_{12}).
\end{aligned}
\end{equation}
Similar to the calculation in Lemma \ref{A.3}, we have that for $I_{11}$
\begin{equation}
\begin{aligned}
    I_{11}\le &C\int_{\Omega^+_1}\frac{\mu^{\frac{N+2}{2}}}{(1+\mu|y-{x}_{1}^{+}|)^{N+2-\alpha} }\frac{k^\alpha}{\mu^\alpha}  |\omega|\\
    \le&C\frac{k^\alpha}{\mu^\alpha}\left(\int_{\Omega^+_1}\frac{\mu^{N}}{(1+\mu|y-{x}_{1}^{+}|)^{2N-\frac{2N\alpha}{N+2}} } \right)^{\frac{N+2}{2N}} \frac{1}{k^{\frac{1}{2^*}}}||\omega||\\
    \le&C \left\{ \begin{matrix}
       \displaystyle \frac{k^\alpha}{\mu^\alpha} \frac{1}{k^{\frac{1}{2^*}}}||\omega||, \qquad\text{if } \ 1<\alpha<\frac{N+2}{2},\\
        \\
        \displaystyle\frac{k^\alpha}{\mu^\alpha}\left(\ln \mu \right)^{\frac{N+2}{2N}} \frac{1}{k^{\frac{1}{2^*}}}||\omega||, \qquad\text{if } \ \alpha = \frac{N+2}{2},
    \end{matrix}\right.\\
    \le&C \frac{1}{k^{\frac{1}{2^*}}}  \left( \frac{k}{\mu} \right)^{\frac{N+2}{2}-\tau} ||\omega||.
\end{aligned}
\end{equation}
For $I_{12},$ we have
\begin{equation}\label{lkI12}
    \begin{aligned}
I_{12}\le&C\int_{\Omega^+_1}\frac{\mu^{\frac{N+2}{2}}}{(1+\mu|y-{x}_{1}^{+}|)^{N+2-\frac{N+2}{N-2}\alpha} }\frac{k^{\frac{N+2}{N-2}\alpha}}{\mu^{\frac{N+2}{N-2}\alpha}}  |\omega|\\
\le& C\left(\frac{k}{\mu}\right) ^{\frac{N+2}{N-2}\alpha}\left(\int_{\Omega^+_1}\frac{\mu^N}{(1+\mu|y-{x}_{1}^{+}|)^{2N-\frac{2N\alpha}{N-2}} }\right)^\frac{N+2}{2N}||\omega||\\
\le&C \frac{1}{k^{\frac{1}{2^*}}}  \left( \frac{k}{\mu} \right)^{\frac{N+2}{2}-\tau} ||\omega||.
    \end{aligned}
\end{equation}

Combining \eqref{lk'}--\eqref{lkI12}, we can obtain \eqref{lk}.
	\end{proof}

Using similar arguments to the proof Lemma 2.4 in \cite{WY}, we have the following estimate of ${\bf N}(\phi_k).$

 \begin{lemma} \label{lem2.3}
		Suppose that $V(|y|)$ satisfies $({\bf V})$, $N\ge 5$ and $(r_k,h_k,\mu_k)\in {\mathscr S_k}$, there exists $C>0$ such that
		\begin{equation}
		    ||{\bf N}(\phi_k)||\le\left\{ \begin{matrix}C\|\phi_k\|^{2^*-1},\quad \text{if} \ N\ge 6, \\

 \\C\displaystyle k^{\frac{1}{10} }\|\phi_k\|^2,\quad \text{if} \ N=5,
\end{matrix}\right.
		\end{equation}
		for all $\phi_k \in \mathbb{E}_k$.
	\end{lemma}

	\medskip
	The  solvability theory  for the  linearized  problem \eqref{lin}  can be provided in the following:
	\begin{proposition} \label{pro2.5}
		Suppose that  $V(|y|)$ satisfies $({\bf V})$, $N\ge 5$ and $(r_k,h_k,\mu_k)\in {\mathscr S_k}$. There exists an integer  $k_0$ large enough, such that for each $k \ge k_0, $ problem \eqref{lin} has a unique solution $\phi_k\in \mathbb{E}_k$ satisfying
		\begin{align} \label{phik}
			\|\phi_k\| \le  Ck\max\left\{ \frac{1}{k^{\frac{1}{2^*}}}  \left( \frac{k}{\mu} \right)^{\frac{N+2}{2}-\tau} , \frac{1}{\mu^{\min\{{\frac{N-2}{2}},{2-\tau}\}}}\right\},
		\end{align}
  where $\tau$ is a small constant. In fact,  $\phi_k$ can be recognized as a $C^1$ map from $ \mathscr {S}_k$ to $\mathbb{E}_k$, which we denote as $\phi_{r,h,\mu}.$ Meanwhile, we have the following estimates for $c_\ell, \ell=1,2,3,$
		\begin{align}\label{cell}
			|c_{\ell}|\le   C \frac{1+\mu^2\delta_{\ell 3}}{\mu(1+r\delta_{\ell 2})}\max\left\{ \frac{k}{k^{\frac{1}{2^*}}}  \left( \frac{k}{\mu} \right)^{\frac{N+2}{2}-\tau} , \frac{k}{\mu^{\min\{{\frac{N-2}{2}},{2-\tau}\}}}\right\}.
		\end{align}
	\end{proposition}
	
	\begin{proof}
		We denote
		\begin{align*}
			\mathcal{A}:= \Bigg\{\omega:  \omega\in \mathbb{E}_k, \quad\|\omega\| \le   C k\max\left\{ \frac{1}{k^{\frac{1}{2^*}}}  \left( \frac{k}{\mu} \right)^{\frac{N+2}{2}-\tau} , \frac{1}{\mu^{\min\{{\frac{N-2}{2}},{2-\tau}\}}}\right\} \Bigg\}.
		\end{align*}
		Then problem \eqref{lin} is equivalent to the following fixed point problem in the sense of distribution:
		\begin{align*}
			\phi_k={\bf A}(\phi_k):= {\bf  L}_k^{-1} \big({{\bf N}}(\phi_k)+{\bf l}_k\big)  ,
		\end{align*}
		where  ${\bf  L}_k^{-1}$ is a linear bounded operator from Lemma \ref{lem2.1}.

		From Lemma \ref{lem2.3} and Lemma \ref{lem2.4}, we know that if $N\ge 6,$
		\begin{align*}
			\|\phi_k\|  & \le C \|{{\bf N}}(\phi_k)\|+C\|{\bf l}_k \|
			\\
			& \le C||\phi_k||^{2^*-1}+Ck\max\left\{ \frac{1}{k^{\frac{1}{2^*}}}  \left( \frac{k}{\mu} \right)^{\frac{N+2}{2}-\tau} , \frac{1}{\mu^{\min\{{\frac{N-2}{2}},{2-\tau}\}}}\right\}
			 \\
			&  \le k\max\left\{ \frac{1}{k^{\frac{1}{2^*}}}  \left( \frac{k}{\mu} \right)^{\frac{N+2}{2}-\tau} , \frac{1}{\mu^{\min\{{\frac{N-2}{2}},{2-\tau}\}}}\right\}.
		\end{align*}
Else if $N=5,$
\begin{align*}
			\|\phi_k\|  & \le C \|{{\bf N}}(\phi_k)\|+C\|{\bf l}_k \|
			\\
			& \le Ck^{\frac{1}{10}}||\phi_k||^{2}+ Ck\max\left\{ \frac{1}{k^{\frac{1}{2^*}}}  \left( \frac{k}{\mu} \right)^{\frac{N+2}{2}-\tau} , \frac{1}{\mu^{\min\{{\frac{N-2}{2}},{2-\tau}\}}}\right\}
			 \\
			&  \le Ck^{\frac{1}{10}}||\phi_k||^{2}+ C\frac{1}{k^{\frac{63}{10}-\frac{2\tau}{N-4}}}\le C\frac{1}{k^{\frac{63}{10}-\frac{2\tau}{N-4}}}.
		\end{align*}
Hence, ${\bf A}$ is a contraction map from  $\mathcal{A}$ to $\mathcal{A}$.

On the other hand, by the similar argument in the proof of Proposition 2.3 in \cite{WY}, for any  $\phi_{k,1}, \phi_{k,2} \in \mathcal{A}$, we can get
	\begin{equation}
 \begin{aligned}
      ||{\bf N}(\phi_{k,1})-{\bf N}(\phi_{k,2})||&\le\left\{ \begin{matrix}C(\|\phi_1\|^{2^*-2}+\|\phi_2\|^{2^*-2})\|\phi_{k,1}-\phi_{k,2}\|,\quad \text{if} \ N\ge 6, \\

 \\C\displaystyle k^{\frac{1}{10} }\|\phi_1-\phi_{k,2}\|^2,\quad \text{if} \ N=5,
\end{matrix}\right.\\
&\le\frac{1}{2}\|\phi_{{k,1}}-\phi_{{k,2}} \|.
 \end{aligned}
		\end{equation}
  Then we get
		\begin{align*}
			\|{\bf A}(\phi_{k,1})-{\bf A}(\phi_{k,2}) \|  \le C\|{{\bf N}}(\phi_{k,1})-{{\bf N}}(\phi_{k2})\|\le\frac{1}{2}\|\phi_{k,1}-\phi_{k,2} \|.
		\end{align*}
        It follows from the contraction mapping principle that there is a unique solution $ \phi_k = {\bf A}(\phi_k)$ in $\mathcal{A}$. Moreover,
        \begin{equation*}
            \|\phi_k\| \le  C||{\bf l}_k||\le Ck\max\left\{ \frac{1}{k^{\frac{1}{2^*}}}  \left( \frac{k}{\mu} \right)^{\frac{N+2}{2}-\tau} , \frac{1}{\mu^{\min\{{\frac{N-2}{2}},{2-\tau}\}}}\right\}.
        \end{equation*}

Next, we give the estimates of $c_\ell,\ell=1,2,3.$
Note that
\begin{equation}
\begin{aligned}
      |\overline{\mathbb Z}_{\ell j}|\le \frac{C(1+\delta_{\ell 2}r)\mu^{\frac{N}{2}}}{\left(1+\mu^2|y-{x}_{j}^{+}|^2\right)^{\frac{N-1}{2}}}\le C (1+\delta_{\ell 2}r)\mu U_{{x}_{j}^{+}, \mu}, \quad\ell=1,2,\qquad|\underline{\mathbb Z}_{3j}|\le\frac{U_{{x}_{j}^{+}, \mu}}{\mu},\\
           |\underline{\mathbb Z}_{\ell j}|\le \frac{C(1+\delta_{\ell 2}r)\mu^{\frac{N}{2}}}{\left(1+\mu^2|y-{x}_{j}^{-}|^2\right)^{\frac{N-1}{2}}}\le C (1+\delta_{\ell 2}r)\mu U_{{x}_{j}^{-}, \mu}, \quad\ell=1,2,\qquad|\underline{\mathbb Z}_{3j}|\le\frac{U_{{x}_{j}^{-}, \mu}}{\mu}.
\end{aligned}
\end{equation}
Calculate the inner product of $\overline{\mathbb Z}_{q 1}, q=1,2,3,$ and the both side of \eqref{lin}, then we get
\begin{equation}\label{xZlj}
\begin{aligned}
    &\int_{\R^N}(-\Delta {\phi_k}+V(|y|)\phi_k-(2^*-1)W_{r,h,\mu}^{2^*-2}\phi_k)\overline{\mathbb Z}_{q 1}\\
    =&\langle{\bf N}(\phi_k)+{\bf l_k},\overline{\mathbb Z}_{q 1}\rangle+\sum\limits_{i=1}^k\sum\limits_{\ell=1}^3 \int_{\R^N}\Big({c}_\ell U_{x^{+}_i,\mu}^{2^*-2}\overline{\mathbb{Z}}_{\ell i}+{c}_\ell U_{x^{-}_i,\mu}^{2^*-2}\underline{\mathbb{Z}}_{\ell i}\Big)\overline{\mathbb Z}_{q 1}.
\end{aligned}
\end{equation}
For the left side of \eqref{xZlj}, we have
\begin{equation*}
    \begin{aligned}
        &\int_{\R^N}(-\Delta {\phi_k}+V(|y|)\phi_k-(2^*-1)W_{r,h,\mu}^{2^*-2}\phi_k)\overline{\mathbb Z}_{q 1}\\
    =&\int_{\R^N}V(|y|)\overline{\mathbb Z}_{q 1}\phi_k-(2^*-1)\int_{\R^N}\left(W_{r,h,\mu}^{2^*-2}-U_{x^{+}_1,\mu}^{2^*-2}\right)\overline{\mathbb Z}_{q 1}\phi_k\\
    :=&I_1-(2^*-1)I_2
    \end{aligned}
\end{equation*}
        Similar to the  calculation in Lemma \ref{lem2.4}, we have
        \begin{equation*}
            \begin{aligned}
I_1=&\int_{B_1(x^{+}_1)}V(|y|)\overline{\mathbb Z}_{q 1}\phi_k+\int_{\R^N\setminus B_1(x^{+}_1)}V(|y|)\overline{\mathbb Z}_{q 1}\phi_k\\
=&\left\{\begin{matrix}
    O\left(\displaystyle\frac{1+r\delta_{q 2}}{\mu^{\min\{\frac{N-4}{2},1-\tau\}}}\right)||\phi_k||, \qquad q=1,2,\\
\\
    O\left(\displaystyle\frac{1}{\mu^{\min\{\frac{N}{2},3-\tau\}}}\right)||\phi_k||, \qquad q=3,
\end{matrix}\right.\\=&O\left(\frac{1+r\delta_{q 2}}{1+\mu^2\delta_{q 3}}\frac{1}{\mu^{\min\{\frac{N-4}{2},1-\tau\}}}\right)||\phi_k||.
            \end{aligned}
        \end{equation*}
And,
\begin{equation*}
    \begin{aligned}
        I_2=\int_{\R^N}\left(W_{r,h,\mu}^{2^*-2}-U_{x^{+}_1,\mu}^{2^*-2}\right)\overline{\mathbb Z}_{q 1}\phi_k=  O\left(k^{\frac{2^*-1}{2^*}}\frac{\mu(1+r\delta_{q 2})}{1+\mu^2\delta_{q 3}}\left(\frac{k}{\mu}\right)^{\frac{N+2-\tau}{2}}\right)||\phi_k||.
    \end{aligned}
\end{equation*}
Thus, \begin{equation}
\begin{aligned}
     &\int_{\R^N}(-\Delta {\phi_k}+V(|y|)\phi_k-(2^*-1)W_{r,h,\mu}^{2^*-2}\phi_k)\overline{\mathbb Z}_{q 1}\\=&\frac{\mu(1+r\delta_{q 2})}{1+\mu^2\delta_{q 3}}O\left(k^{\frac{2^*-1}{2^*}}\left(\frac{k}{\mu}\right)^{\frac{N+2-\tau}{2}}+\frac{1}{\mu^{\min\{\frac{N-2}{2},2-\tau\}}}\right)||\phi_k||.
\end{aligned}
\end{equation}
On the other hand,
\begin{equation}
    \begin{aligned}
        \langle{\bf N}(\phi_k)+{\bf l_k},\overline{\mathbb Z}_{q 1}\rangle =
       \frac{\mu(1+r\delta_{q 2})}{1+\mu^2\delta_{q 3}} O\left(\frac{1}{\mu^{\min\{\frac{N-2}{2},2-\tau\}}}\right)||{\bf N}(\phi_k)+{\bf l_k}||.
    \end{aligned}
\end{equation}
Note that
\begin{equation*}
    \begin{aligned}
        \int_{\R^N} U_{x^{+}_1,\mu}^{2^*-2}\overline{\mathbb{Z}}_{\ell 1}\overline{\mathbb Z}_{q 1}=\left\{\begin{matrix}
            0, \qquad\qquad \qquad \ \ell\neq q,\\
\\
\bar {c_\ell}\displaystyle\frac{\mu^2(1+r^2\delta_{\ell 2})}{1+\mu^4\delta_{\ell 3}},\ \ \ell=q.
        \end{matrix}\right.
    \end{aligned}
\end{equation*}
Then,
\begin{equation}\label{ellq}
    \begin{aligned}
        \sum\limits_{i=1}^k\sum\limits_{\ell=1}^3 \int_{\R^N}\Big({c}_\ell U_{x^{+}_i,\mu}^{2^*-2}\overline{\mathbb{Z}}_{\ell i}+{c}_\ell U_{x^{-}_i,\mu}^{2^*-2}\underline{\mathbb{Z}}_{\ell i}\Big)\overline{\mathbb Z}_{q 1}={c_q}\bar{c_q}\displaystyle\frac{\mu^2(1+r^2\delta_{q 2})}{1+\mu^4\delta_{q 3}}(1+o(1)),
    \end{aligned}
\end{equation}
for some constant $\bar {c_q}>0,q=1,2,3.$

       Therefore, combining \eqref{xZlj}--\eqref{ellq}, we obtain
       \begin{equation}\label{cq}
       \begin{aligned}
           |c_{q}|=&\frac{1+\mu^2\delta_{q 3}}{\mu(1+r\delta_{q 2})}O\left(k^{\frac{2^*-1}{2^*}}\left(\frac{k}{\mu}\right)^{\frac{N+2-\tau}{2}}+\frac{1}{\mu^{\min\{\frac{N-2}{2},2-\tau\}}}\right)||\phi_k||\\&+\frac{1+\mu^2\delta_{q 3}}{\mu(1+r\delta_{q 2})} O\left(\frac{1}{\mu^{\min\{\frac{N-2}{2},2-\tau\}}}\right)||{\bf N}(\phi_k)+{\bf l_k}||\\
           \le& C \frac{1+\mu^2\delta_{q 3}}{\mu(1+r\delta_{q 2})}||{\bf l_k}||,
       \end{aligned}
       \end{equation}
       then we have proved \eqref{cell}.
	\end{proof}

	\medskip

\medskip
	\section{Proof of the Main Results}  \label{sec3}

	\begin{proposition} \label{pro3.1}
		Let $ \phi_{r,h, \mu}$ be a function obtained in Proposition \ref{pro2.5} and
		\begin{align*}
			F(r,h,\mu):= I( W_{r, h, \mu}+\phi_{r,h, \mu}).
		\end{align*}

		If $(r,h, \mu)$ is a critical point of  $F(r,h,\mu)$, then
		\begin{align*}
			u= W_{r,h,\mu}+\phi_{r,h,\mu}
		\end{align*}
		is a critical point of $I(u)$ in  $\overline{H^1(\R^N) \cap D^{1,2}(\R^N)}$. \qed
	\end{proposition}
	
	\medskip
	First, we will give the expression of  $F(r,h,\mu)$. 
	
	\begin{proposition}\label{pro3.2}
		Suppose that  $V(|y|)$ satisfies $({\bf V})$, $N\ge 5$ and $(r_k,h_k,\mu_k)\in {\mathscr S_k}$.  We have the following expansion as $k \to +\infty$:
 \begin{equation}\label{F(r,h,mu)}
\begin{aligned}
F(r,h,\mu)=&I(W_{r, h, \mu})+ko\left(\frac{1}{\mu^2}\right) \\=&k\left(A_1+\frac{A_2 V(r)}{\mu^2}-\frac{B_3 k^{N-2}}{\mu^{N-2}(\sqrt{1-h^2})^{N-2}}- \frac{B_4 k }{\mu^{N-2}h^{N-3}\sqrt{1-h^2}} \right)\\
    &+k o\left(\frac{1}{\mu^2}\right)+k O\left(\frac{k^{N-\sigma}}{\mu^{N-\sigma}}+\frac{k^{N-4} \ln k}{\mu^{N-2}}\right)+k O\left(\frac{1}{\mu^{N-2}h^{N-2}}\right),
\end{aligned}
\end{equation}
where $B_3, B_4$ are positive constants defined in Proposition \ref{A.4} and  $ \sigma $ is a  constant small enough.
	\end{proposition}
	\begin{proof}
        By direct calculation, we first have
        \begin{equation}\label{F(1)}
            \begin{aligned}
                F(r,h,\mu)=&I(W_{r, h, \mu})+\frac{1}{2}\int_{\R^N} (|D\phi_{r,h,\mu}|^2+V(|y|)\phi_{r,h,\mu}^2)\\
                &-\int_{\R^N} \left(W_{r, h, \mu}^{2^*-1}-\sum_{j=1}^k \left(U_{{x}_{j}^{+}, \mu}^{2^*-1}+U_{{x}_{j}^{-}, \mu}^{2^*-1}\right)-V(|y|)W_{r, h, \mu} \right)\phi_{r,h,\mu}\\
                &-\frac{1}{2^*}\int_{\R^N}\left( \left(W_{r, h, \mu}+\phi_{r,h,\mu}\right)^{2^*}-W_{r, h, \mu}^{2^*}-2^*W_{r, h, \mu}^{2^*-1}\phi_{r,h,\mu}\right)\\
                =&I(W_{r, h, \mu})+\frac{1}{2}\langle{\bf L_k}\phi_{r, h, \mu},\phi_{r, h, \mu} \rangle-\langle{\bf l_k},\phi_{r, h, \mu}\rangle-\frac{1}{2^*}{\bf R}(\phi_{r, h, \mu})\\
                \le&I(W_{r, h, \mu})+C\left(||\phi_{r, h, \mu}||^2+||{\bf l_k}|| ||\phi_{r, h, \mu}|| +||{{\bf N}(\phi_{r,h,\mu})}|| ||W_{r, h, \mu}+\phi_{r,h,\mu}||
                \right)\\&+C\left(\int_{\R^N}W_{r, h, \mu}^{2^*}\right)^{\frac{2}{N}}||\phi_{r,h,\mu}||^2\\
                =&I(W_{r, h, \mu})+O\left(k^{\frac{2}{N}}||\phi_{r, h, \mu}||^2\right)+O\left(k^{{\frac{1}{2^*}}}||\phi_{r,h,\mu}||^{2^*}+k^{\frac{13}{10}}||\phi_{r,h,\mu}||^{3}\right)\\
                =&I(W_{r, h, \mu})+O\left(k^{\frac{N+4}{N}}\left(\frac{k}{\mu}\right)^{N+2-2\tau}+k^{\frac{2N+2}{N}}\frac{1}{\mu^{\min\{N-2,4-2\tau\}}}\right),
            \end{aligned}
        \end{equation}
where
\begin{equation}
\begin{aligned}
       {\bf R}(\phi_{r, h, \mu}):=&\int_{\R^N}\left( \left(W_{r, h, \mu}+\phi_{r,h,\mu}\right)^{2^*}-W_{r, h, \mu}^{2^*}-2^*W_{r, h, \mu}^{2^*-1}\phi_{r,h,\mu}-\frac{2^*(2^*-1)}{2}W_{r, h, \mu}^{2^*-2}\phi_{r,h,\mu}^2\right)\\
       =&\langle {{\bf N}(\phi_{r,h,\mu})},W_{r, h, \mu}+\phi_{r,h,\mu}\rangle-\frac{2(2^*-1)}{N-2}\int_{\R^N}W_{r, h, \mu}^{2^*-2}\phi_{r,h,\mu}^2.
\end{aligned}
\end{equation}
Note that
\begin{equation}
    \begin{aligned}
        ||\phi_{r,h,\mu}||\le Ck\left\{\begin{matrix}
        \displaystyle\frac{1}{k^{\frac{3N^2-2N+8}{2N(N-4)}}}\qquad N\ge 8,\\
        \\
        \displaystyle\frac{1}{k^{\frac{4N^2-8
N}{2N(N-4)}}},\qquad N=6,7,\\
        \\
           \displaystyle \frac{1}{k^\frac{9}{2}},\qquad N=5.
        \end{matrix}\right.
    \end{aligned}
\end{equation}
Then if $N \ge 8,$
\begin{equation}\label{o1}
    k^{\frac{2}{N}}||\phi_{r, h, \mu}||^2=kO\left(\frac{1}{k^{\frac{2N^2-16}{N(N-4)}}}\right)=ko\left(\frac{1}{k^{\frac{2N^2-4N}{N(N-4)}}}\right)=ko\left(\frac{1}{\mu^2}\right).
\end{equation}
Else if $N=5,6,7$
\begin{equation}\label{o2}
\begin{aligned}
       k^{\frac{2}{N}}||\phi_{r, h, \mu}||^2=k\left\{\begin{matrix}
O\left(k^{-\frac{113}{21}}\right),\qquad N=7,\\
\\
O\left(k^{-\frac{20}{3}}\right),\qquad N=6,\\
\\
O\left(k^{-\frac{38}{5}}\right),\qquad N=5,
       \end{matrix}
       \right.
       =k\left\{\begin{matrix}
o\left(k^{-\frac{10}{3}}\right),\qquad N=7,\\
\\
o\left(k^{-4}\right),\qquad N=6,\\
\\
o\left(k^{-6}\right),\qquad N=5,
       \end{matrix}
       \right.
    =ko\left(\frac{1}{\mu^2}\right).
\end{aligned}
\end{equation}

Combining \eqref{F(1)},\eqref{o1},\eqref{o2} and Proposition \ref{A.4}, we finally get \eqref{F(r,h,mu)}.

		\end{proof}

		\medskip
		Next, we will give the expansions of $\displaystyle\frac{\partial F(r,h, \mu)}{\partial h} $ and $\displaystyle\frac{\partial F(r,h, \mu)}{\partial \mu} $.
		\begin{proposition}\label{pro3.3}
			Suppose that  $V(|y|)$ satisfies $({\bf V})$, $N\ge 5$ and $(r_k,h_k,\mu_k)\in {\mathscr S_k}$.  We have the following expansion as $k \to +\infty$:
   \begin{equation}\label{Fh(W)}
\begin{aligned}
    \frac{\partial F(r,h, \mu)}{\partial h}=&k\left(-\frac{(N-2)B_3 k^{N-2}h}{\mu^{N-2}(\sqrt{1-h^2})^{N}}+\frac{(N-3)B_4 k}{\mu^{N-2}h^{N-2}\sqrt{1-h^2}}-\frac{B_4 k}{\mu^{N-2}h^{N-4} (\sqrt{1-h^2})^3 }\right)
\\&+kO\left(  \frac{k^{N-3}\ln k}{\mu^{N-2}}\right)+kO\left(  \frac{1}{\mu^{N-2}h^{N-1}}+\frac{hk^{N-\sigma}}{\mu^{N-\sigma}}\right),
\end{aligned}
\end{equation}
and
\begin{equation}\label{Fmu(W)}
\begin{aligned}
    \frac{\partial F(r,h, \mu)}{\partial \mu}=&k\left(-\frac{2 A_2 V(r)}{\mu^3}+\frac{(N-2)B_3k^{N-2}}{\mu^{N-1}(\sqrt{1-h^2})^{N-2}}+ \frac{(N-2) B_4 k}{\mu^{N-1}h^{N-3}\sqrt{1-h^2}} \right)\\
    &+k O\left(\frac{k^{N-\sigma}}{\mu^{N+1-\sigma}}+\frac{k^{N-3}\ln k}{\mu^{N-1}}\right)+k O\left(\frac{1}{\mu^{N-1}h^{N-2}}\right),
\end{aligned}
\end{equation}
where $B_3, B_4$ are positive constants in Lemma \ref{A.4} and  $ \sigma $ is a  constant small enough.
		
		\end{proposition}
		\begin{proof}
			We have
   \begin{equation}\label{Fh11}
       	\begin{aligned}
				\frac{\partial F(r,h, \mu)}{\partial h }
				&    = \left\langle  I'( W_{r, h, \mu}+\phi_{r,h, \mu}), \frac{\partial  W_{r, h, \mu}  } {\partial  h} \right\rangle
				+  \left\langle  I'( W_{r, h, \mu}+\phi_{r,h, \mu}), \frac{\partial  \phi_{r,h, \mu}  } {\partial  h}\right\rangle
				\\
				&    =   \left\langle  I'( W_{r, h, \mu}+\phi_{r,h, \mu}), \frac{\partial  W_{r, h, \mu}  } {\partial  h} \right\rangle  +\sum\limits_{j=1}^k\sum\limits_{\ell=1}^3 {c_\ell}\left\langle
				 \overline{\mathbb{Z}}_{\ell j}+\underline{\mathbb{Z}}_{\ell j}  ,\frac{\partial  \phi_{r,h, \mu}  } {\partial  h}\right\rangle.	\end{aligned}
   \end{equation}
			Noting that $ \displaystyle{\int_{{\mathbb{R}}^N} U_{x^{+}_{j}, \mu}^{2^*-2}  \overline{\mathbb{Z}}_{\ell j}  \phi_{r,h, \mu}}=  \displaystyle{\int_{{\mathbb{R}}^N} U_{x^{-}_{j}, \mu}^{2^*-2}  \underline{\mathbb{Z}}_{\ell j}  \phi_{r,h, \mu}=0}$, then
			\begin{align*}
		\left\langle\overline{\mathbb{Z}}_{\ell j}, \frac{\partial  \phi_{r,h, \mu}  } {\partial  h} \right\rangle=-\left\langle\frac{\partial{\overline{\mathbb{Z}}_{\ell j}}}{\partial h},\phi_{r,h, \mu}\right\rangle,\qquad \left\langle\underline{\mathbb{Z}}_{\ell j}, \frac{\partial  \phi_{r,h, \mu}  } {\partial  h} \right\rangle=-\left\langle\frac{\partial{\underline{\mathbb{Z}}_{\ell j}}}{\partial h},\phi_{r,h, \mu}\right\rangle.
			\end{align*}
			Thus, using Proposition \ref{pro2.5}, we have

   \begin{equation}
       \begin{aligned}
    \left|\sum\limits_{j=1}^k\sum\limits_{\ell=1}^3 {c_\ell}\left\langle
				 \overline{\mathbb{Z}}_{\ell j}+\underline{\mathbb{Z}}_{\ell j}  ,\frac{\partial  \phi_{r,h, \mu}  } {\partial  h}\right\rangle\right|
    \le&\sum\limits_{\ell=1}^3 |{c_\ell}| \left(\left\|\frac{\partial{\overline{\mathbb{Z}}_{\ell j}}}{\partial h}\right\|+\left\|\frac{\partial{\underline{\mathbb{Z}}_{\ell j}}}{\partial h}\right\|\right)||\phi_{r,h, \mu}||\\[2mm]
    \le&Cr \mu ||{\bf l_k} || ||\phi_{r,h, \mu}||=\mu O\left(||\phi_{r,h, \mu}||^2 \right)\\[2mm]
    =&o\left(\frac{\ln k}{\mu^3}\right)=ko\left(\frac{k^{N-3}\ln k}{\mu^{N-2}}\right).
       \end{aligned}
   \end{equation}
		
			On the other hand, similar to the calculation in Lemma \ref{lem2.4}, we have
   \begin{equation}\label{Fh12}
       \begin{aligned}
& \left\langle  I'( W_{r, h, \mu}+\phi_{r,h, \mu}), \frac{\partial  W_{r, h, \mu}  } {\partial  h} \right\rangle\\
=&\left\langle  I'( W_{r, h, \mu}), \frac{\partial  W_{r, h, \mu}  } {\partial  h}\right\rangle+\int_{\R^N}V(|y|)\frac{\partial  W_{r, h, \mu}  } {\partial  h}\phi_{r,h, \mu}-\int_{\R^N}\left(( W_{r, h, \mu}+\phi_{r,h, \mu})^{2^*-1}-W_{r, h, \mu}^{2^*-1}\right)\frac{\partial  W_{r, h, \mu}  } {\partial  h}\\
=&\left\langle  I'( W_{r, h, \mu}), \frac{\partial  W_{r, h, \mu}  } {\partial  h}\right\rangle+r\mu O( ||{\bf l_k}||||\phi_{r,h, \mu}||)+O\left(||\phi_{r,h, \mu}||^{2^*}\right)\\
=&\left\langle  I'( W_{r, h, \mu}), \frac{\partial  W_{r, h, \mu}  } {\partial  h}\right\rangle+ko\left(\frac{k^{N-3}\ln k}{\mu^{N-2}}\right).
       \end{aligned}
   \end{equation}
Combining \eqref{Fh11}--\eqref{Fh12} and Proposition \ref{A.6}, we finally get \eqref{Fh(W)}.

Similarly, we have
\begin{equation}\label{Fmu11}
    \begin{aligned}
        \frac{\partial F(r,h, \mu)}{\partial \mu }
				    =& \left\langle  I'( W_{r, h, \mu}+\phi_{r,h, \mu}), \frac{\partial  W_{r, h, \mu}  } {\partial  \mu} \right\rangle
				+  \left\langle  I'( W_{r, h, \mu}+\phi_{r,h, \mu}), \frac{\partial  \phi_{r,h, \mu}  } {\partial  \mu}\right\rangle
				\\
				=& \left\langle  I'( W_{r, h, \mu}), \frac{\partial  W_{r, h, \mu}  } {\partial  \mu}\right\rangle+\frac{1}{\mu} O( ||{\bf l_k}||||\phi_{r,h, \mu}||)+O\left(||\phi_{r,h, \mu}||^{2^*}\right)\\
    =& \left\langle  I'( W_{r, h, \mu}), \frac{\partial  W_{r, h, \mu}  } {\partial  \mu}\right\rangle+ko\left(\frac{k^{N-3}\ln k}{\mu^{N-1}}\right).
    \end{aligned}
\end{equation}
Then combining \eqref{Fmu11} and Proposition \ref{A.6}, we obtain \eqref{Fmu(W)}.
\end{proof}

		\medskip
Let $\mu_0$ be the solution of
  \[-\frac{2 A_2 V(r)}{\mu^3}+\frac{(N-2)B_3}{\mu^{N-1}}=0,\]
  which gives that
		\begin{align} \label{mu0}
			\mu_0(r)=\left(\frac{(N-2) B_0B_1}{2A_2 V(r)r^{N-2}}\right
   )^{\frac 1 {N-4}}.
		\end{align}
		And let $h_0$ be the solution of
  $$-{(N-2)B_3 h}+\frac{(N-3)B_4 }{h^{N-2}}=0,$$
  which gives that
		\begin{align} \label{h0}
			h_0=\left(\frac{(N-3) B_2}{(N-2) B_1}\right
   )^{\frac 1 {N-1}}.
		\end{align}
  \medskip

 Denote
\begin{equation}
    \mathcal{F}( h,\mu):=\frac{A_2 V(r)}{\mu^2}-\frac{B_3 k^{N-2}}{\mu^{N-2}(\sqrt{1-h^2})^{N-2}}- \frac{B_4 k }{\mu^{N-2}h^{N-3}\sqrt{1-h^2}}.
\end{equation}
Expanse $\mathcal{F}(h,\mu)$ at the point $(H_0, \Lambda_0(r)),$ where $H_0:=\displaystyle\frac{h_0}{k^{\frac{N-3}{N-1}}},\Lambda_0(r):=\mu_0(r)k^{\frac{N-2}{N-4}},$ we have
\begin{equation}\label{Ff}
    \begin{aligned}
    \mathcal{F}( h,\mu)
    =&\mathcal{F}( H_0,\Lambda_0)+O\left(\frac{k}{\Lambda_0^{N-2}H_0^{N-4}}\right)|h-H_0|+O\left(\frac{k}{\Lambda_0^{N-1}H_0^{N-3}}\right)|\mu-\Lambda_0|\\&+O\left(\frac{k}{\Lambda_0^{N-2}H_0^{N-1}}\right)|h-H_0|^2+O\left(\frac{k^{N-2}}{\Lambda_0^N}\right)|\mu-\Lambda_0|^2+O\left(\frac{k}{\Lambda_0^{N-1}H_0^{N-1}}\right)|h-H_0||\mu-\Lambda_0|,\\
\end{aligned}
\end{equation}
  where $(H_1,\Lambda_1)$ are the mid-value between $(h,\mu)$ and $(H_0,\Lambda_0).$\\ \\

According to \eqref{Ff}, we denote $\tilde \sigma=\frac{1}{k^\frac{2(N^2-7N+9)}{(N-4)(N-1)}},$ and define
 \begin{align}\label{Skk}
			{{\mathbb S}_k}
			= \Bigg\{&(r,h,\mu) :\,  r\in \Big[r_0-\tilde \sigma, r_0+\tilde  \sigma \Big], \quad
			h \in \Bigg[\frac{h_0-\tilde \sigma}{k^{\frac{N-3}{N-1}}}  ,
			\frac{h_0+\tilde \sigma}{k^{\frac{N-3}{N-1}}} \Bigg], \nonumber
			\\[2mm]
			& \qquad \qquad
			\mu \in \Big[(\mu_0-\tilde\sigma) k^{\frac{N-2}{N-4}},  (\mu_0+\tilde \sigma) k^{\frac{N-2}{N-4}}\Big],\quad \tilde \sigma=\frac{1}{k^\frac{2(N^2-7N+9)}{(N-4)(N-1)}}\Bigg\}.
		\end{align}
Next, we will find a critical point to $F(r,h,\mu)$ in the interior of ${{\mathbb S}_k}$.\\ \medskip
\\
{\bf{Proof of Theorem \ref{main1}} if ${\bf r_0}$ is a maximum of ${\bf r^2V(r)}$:} Noting that
\begin{align}
		&\frac{\partial F(r,h,\mu)}{\partial h}\bigg(r,\frac{h_0-\tilde \sigma}{k^{\frac{N-3}{N-1}}}, \mu\bigg)>0,\quad \frac{\partial F(r,h,\mu)}{\partial h}\bigg(r, \frac{h_0+\tilde \sigma}{k^{\frac{N-3}{N-1}}}, \mu\bigg)<0, \quad \nonumber\\[2mm] \text{for any}\ (r,\mu) \in&\left[r_0-\tilde \sigma, r_0+\tilde  \sigma \right]\times\left[(\mu_0-\tilde\sigma) k^{\frac{N-2}{N-4}},  (\mu_0+\tilde \sigma) k^{\frac{N-2}{N-4}} \right],  \nonumber\\[2mm]
		&\frac{\partial F(r,h,\mu)}{\partial \mu}\bigg(r, h,(\mu_0-\tilde\sigma) k^{\frac{N-2}{N-4}}\bigg)>0,\quad \frac{\partial F(r,h,\mu)}{\partial \mu}\bigg(r, h,(\mu_0+\tilde\sigma) k^{\frac{N-2}{N-4}}\bigg)<0,\nonumber\\ \text{for any}\ (r,h) \in&\left[r_0-\tilde \sigma, r_0+\tilde  \sigma \right]\times\left[\frac{h_0-\tilde \sigma}{k^{\frac{N-3}{N-1}}}  ,
			\frac{h_0+\tilde \sigma}{k^{\frac{N-3}{N-1}}}  \right].
	\end{align}
On the other hand, for any $(r,h,\mu)\in\mathbb{S}_k,$ using \eqref{Ff} we have
\begin{equation}\label{FFF}
    \begin{aligned}
        F(r,h,\mu)=&k\left(A_1+\frac{A_2 V(r)}{\mu^2}-\frac{B_3 k^{N-2}}{\mu^{N-2}(\sqrt{1-h^2})^{N-2}}- \frac{B_4 k }{\mu^{N-2}h^{N-3}\sqrt{1-h^2}} \right)\\
    &+k o\left(\frac{1}{\mu^2}\right)+k O\left(\frac{k^{N-\sigma}}{\mu^{N-\sigma}}+\frac{k^{N-4} \ln k}{\mu^{N-2}}\right)+k O\left(\frac{1}{\mu^{N-2}h^{N-2}}\right)\\
   =&k\left(\left(A_1+A_3 (r^2 V(r))^{\frac{N-2}{2(N-4)}}+o(1)\right)\frac{1}{k^{\frac{2(N-2)}{N-4}}}+O\left(\frac{1}{k^\frac{N^2-7N+9}{(N-4)(N-1)}}+o(1) \right)\right)\\
   =&k\left(A_1+A_3 (r^2 V(r))^{\frac{N-2}{2(N-4)}}+o(1)\right)\frac{1}{k^{\frac{2(N-2)}{N-4}}},\\
    \end{aligned}
\end{equation}
where $A_3:=\displaystyle\frac{(N-4)A_2}{N-2}\left(\frac{2A_2}{(N-2)B_0B_1}\right)^\frac{2}{N-4}.$\\

Since $r^2v(r)$ has a maximum at $r_0$, form \eqref{FFF}, we deduce that if $r=r_0 \pm \tilde \sigma, F(r,h,\mu)$ can not reaches its maximum. That is, the maximum point of $F(r,h,\mu)$, denoted as $(\bar {r_k},\bar{h_k},\bar{\mu_k}),$  is an interior point of $\mathbb{S}_k$. Thus, $(\bar {r_k},\bar{h_k},\bar{\mu_k})$ is a critical point of $F(r,h,\mu)$.
\qed
\medskip

  Before proving the minimum case, we need to establish a gradient flow system and study on its flow. Let
		\begin{align}\label{barF}
			\bar  F(r,h,\mu)=-F(r,h,\mu),\qquad (r,h,\mu)\in \mathbb{S}_k.
		\end{align}
		and
		\[t_2= k (-A_1+\eta_0),
		\quad t_1
		= k  \left(-A_1 -A_3 (1+\eta)(r_0^2 V(r_0)) ^{\frac{N-2}{2(N-4)}}\frac{1}{k^{\frac{2(N-2)}{N-4}}}\right),  \]
		where $\eta_0>0 $ is a small constant and $\eta>0$  satisfying $V(r_0\pm \tilde\sigma)>(1+\eta)V(r_0)$. We also define the energy level set
		\[  \bar F^{t}= \Big\{(r,h,\mu) \big|~(r,h,\mu) \in {D_k}, ~ \bar F (r,h,\mu) \le{t} \Big\}.  \]
		
		Consider the following  gradient flow system
		\begin{equation*}
			\begin{cases}
				\displaystyle\frac{{\mathrm d} r}{{\mathrm d} t}=-{\bar F}_r , &t>0;\\[4mm]
				\displaystyle\frac{{\mathrm d} h}{{\mathrm d} t}=-{\bar F}_h , &t>0;\\[4mm]
				\displaystyle\frac{{\mathrm d} \mu}{{\mathrm d} t}=-{\bar F}_\mu, &t>0;\\[4mm]
				(r,h,\mu) \big|_{t=0}  \in \bar F^{t_2}.
			\end{cases}
		\end{equation*}
		Then we have the following proposition:
		
		\begin{proposition} \label{pro3.4}
			The flow would not leave  $ \mathbb{S}_k$ before it reaches  $ \bar F^{t_1}. $
		\end{proposition}
		\begin{proof}
			There are   three  cases that the flow tends to leave  $ \mathbb{S}_k$:
			\\[2mm]
			\item[ {\bf Case  1.} ]
			$\left|h-\displaystyle\frac{h_0}{k^{\frac{N-3}{N-1}}}\right|= \displaystyle\frac{\tilde \sigma}{k^{\frac{N-3}{N-1}}} $ and  $|r-r_0|\le\tilde\sigma,\quad|\mu-{\mu_0}k^{\frac{N-2}{N-4}}| \le \tilde\sigma k^{\frac{N-2}{N-4}} $;
			\\[2mm]
			~\item[ {\bf Case  2.} ] $|\mu-{\mu_0}k^{\frac{N-2}{N-4}}|= \tilde\sigma k^{\frac{N-2}{N-4}}   $  when  $|r-r_0|\le\tilde\sigma, \quad\left|h-\displaystyle\frac{h_0}{k^{\frac{N-3}{N-1}}}\right| \le\displaystyle\frac{\tilde \sigma}{k^{\frac{N-3}{N-1}}}$;
			\\[2mm]
			~ \item[ {\bf Case  3.} ]  $|r-r_0|=\tilde\sigma $ when  $\left|h-\displaystyle\frac{h_0}{k^{\frac{N-3}{N-1}}}\right| \le\displaystyle\frac{\tilde \sigma}{k^{\frac{N-3}{N-1}}}, \quad|\mu-{\mu_0}k^{\frac{N-2}{N-4}}| \le \tilde\sigma k^{\frac{N-2}{N-4}} . $
			\\[2mm]

			First consider {\bf Case 1} . Noting that
   \begin{equation*}
       \frac{\partial {\bar F}(r,h,\mu)}{\partial h}\bigg(r,\frac{h_0-\tilde \sigma}{k^{\frac{N-3}{N-1}}}, \mu\bigg)<0,\quad \frac{\partial {\bar F}(r,h,\mu)}{\partial h}\bigg(r, \frac{h_0+\tilde \sigma}{k^{\frac{N-3}{N-1}}}, \mu\bigg)>0,
   \end{equation*}
for any $ (r,\mu) \in\left[r_0-\tilde \sigma, r_0+\tilde  \sigma \right]\times\left[(\mu_0-\tilde\sigma) k^{\frac{N-2}{N-4}},  (\mu_0+\tilde \sigma) k^{\frac{N-2}{N-4}} \right],$ then the flow can not leave $\mathbb{S}_k$.

\medskip
			In {\bf Case 2}, noting that
			\begin{equation*}
			    \frac{\partial {\bar F}(r,h,\mu)}{\partial \mu}\bigg(r, h,(\mu_0-\tilde\sigma) k^{\frac{N-2}{N-4}}\bigg)<0,\quad \frac{\partial {\bar F}(r,h,\mu)}{\partial \mu}\bigg(r, h,(\mu_0+\tilde\sigma) k^{\frac{N-2}{N-4}}\bigg)>0,
			\end{equation*}
				for any $ (r,h) \in \left[r_0-\tilde \sigma, r_0+\tilde  \sigma \right]\times\left[\displaystyle\frac{h_0-\tilde \sigma}{k^{\frac{N-3}{N-1}}}  ,
			\frac{h_0+\tilde \sigma}{k^{\frac{N-3}{N-1}}}  \right],$ then the flow can not leave $\mathbb{S}_k$ neither.
   \medskip

    Finally, we consider {\bf Case 3}.
				If  $r=r_0\pm\tilde\sigma $, then from
				\eqref{FFF} and the assumptions that $V(|y|)$ satisfying $({\bf V})$, we can obtain
				\begin{equation}\label{t1}
				    \begin{aligned}
			\bar F(r,h,\mu)=&k\left(-A_1-A_3 ((r_0\pm\tilde\sigma)^2 V(r_0\pm\tilde\sigma))^{\frac{N-2}{2(N-4)}}+o(1)\right)\frac{1}{k^{\frac{2(N-2)}{N-4}}}+kO\left(\frac{1}{k^\frac{2(N^2-7N+9)}{(N-4)(N-1)}}\right)\\
   =&k\left(-A_1-A_3 ( (1+\eta)(r_0^2V(r_0))^{\frac{N-2}{2(N-4)}}+O(r_0\tilde\sigma)+o(1)\right)\frac{1}{k^{\frac{2(N-2)}{N-4}}}+kO\left(\frac{1}{k^\frac{N^2-7N+9}{(N-4)(N-1)}}\right)\\
   =&k\left(-A_1-A_3 ((1+\eta)(r_0^2V(r_0))^{\frac{N-2}{2(N-4)}}+o(1)\right)\frac{1}{k^{\frac{2(N-2)}{N-4}}}+kO\left(\frac{1}{k^\frac{N^2-7N+9}{(N-4)(N-1)}}+o(1)\right)\\
   <&
		 k  \Big(-A_1 -A_3 (1+\eta)(r_0^2 V(r_0)) ^{\frac{N-2}{2(N-4)}}\frac{1}{k^{\frac{2(N-2)}{N-4}}}\Big)=t_1,
				    \end{aligned}
				\end{equation}
				with $\eta>0$ small enough.
				
				Combining above results, we conclude that the flow would not leave  $\mathbb{S}_k$ before it reach  $ {\bar{F}^{t_1}}$.
			\end{proof}	
   \medskip
   Proposition \ref{pro3.4} plays an important role in the proof Theorem \ref{main1}, now let us give the proof of it.\\
   {\bf{Proof of Theorem \ref{main1}} if ${\bf r_0}$ is a minimum of ${\bf r^2V(r)}$:} Define
			\[
			\begin{split}
				G=\Bigl\{\gamma :\quad &\gamma(r,h,\mu)= \big(\gamma _1(r,h,\mu),\gamma _2(r,h,\mu),\gamma _3(r,h,\mu)\big)\in{\mathbb{S}_k},(r,h,\mu)\in{\mathbb{S}_k}; \\[2mm]
				&\gamma(r,h,\mu)=(r,h,\mu), \;\text{if} ~|r-r_0|=  \tilde\sigma \Bigr\}.
			\end{split}
			\]
			
			Let
			\[
			{\bf c}=\inf_{\gamma  \in G}\max_{(r,h,\mu)\in{\mathbb{S}_k}} \bar F \big(\gamma(r,h,\mu)\big).
			\]
			
			We claim that  $ {\bf c}  $ is a critical value of  $ \bar F(r,h, \mu)  $ and can be achieved by some  $(r,h, \mu) \in{\mathbb{S}_k} $. By the minimax theory, we need to show that
			\medskip
			\begin{itemize}
				
				\item[(i)]  $ t_1< {\bf c} < t_2 $;
				\medskip

				\item[(ii)] \label{itemii}  $\underset{{|r-r_0|=  \tilde\sigma}}{\sup}
				{\bar F}\big(\gamma(r,h,\mu)\big)<t_1,\;\forall\; \gamma \in G. $
				
			\end{itemize}
			
			To prove (ii),  let $\gamma\in G$. Then for any  $\bar r$ with
			$|\bar r-r_0|=  \tilde\sigma$, we have $\gamma(\bar r,h,\mu)= (\bar r,h,\mu)$.
			Thus, by \eqref{t1},
			
			\[
			\bar F(\gamma(\bar r,h,\mu))= \bar F(\bar  r,h, \mu)<t_1.
			\]
			
			Now we prove (i).  It is obvious that
			$
			{\bf c}<t_2.
			$
			
			For any $\gamma=(\gamma_1,\gamma_2,\gamma_3)\in G$, we have $\gamma_1(r,h,\mu)= r$, if
			$| r-r_0|=  \tilde\sigma$. Define
			
			\[
			\tilde \gamma_1(r)=\gamma_1\left(r,\displaystyle\frac{h_0}{k^{\frac{N-3}{N-1}}},\mu_0k^{\frac{N-2}{N-4}}\right).
			\]

			Then  $\tilde \gamma_1(r)= r$, if $| r-r_0|=  \tilde\sigma$.
			So, there is a $\bar r\in (r_0-\tilde\sigma,r_0+\tilde\sigma)$, such that
			$
			\tilde \gamma_1(\bar r)={r_0}.
			$

			Let $\bar h= \gamma_2\left(r,\displaystyle\frac{h_0}{k^{\frac{N-3}{N-1}}},\mu_0k^{\frac{N-2}{N-4}}\right),\bar \mu=\gamma_3\left(r,\displaystyle\frac{h_0}{k^{\frac{N-3}{N-1}}},\mu_0k^{\frac{N-2}{N-4}}\right)$. Then from \eqref{F(r,h,mu)} and \eqref{FFF}
			\begin{align*}
				&\max_{(r,h,\mu)\in {\mathbb{S}_k}}\bar F(\gamma(r,h,\mu))\ge \bar F(\gamma(\bar r,h_0,\mu_0))=
				\bar F(r_0, \bar h,\bar \mu)\\
				=& k\left(-A_1-A_3 (r_0^2 V(r_0))^{\frac{N-2}{2(N-4)}}+o(1)\right)\frac{1}{k^{\frac{2(N-2)}{N-4}}}>t_1.
			\end{align*}
\qed

   \medskip

\appendix

\section{Basic estimates and lemmas}\label{secB}
\begin{lemma}\label{A.1}
   Assume that $\displaystyle\frac{k}{\mu\sqrt{1-h^2}}=o(1), \mu h \to +\infty,N\ge 5$, then we have the following  expansion, for $k \to +\infty$:
\begin{equation}\label{2*-1++}
    \int_{\R^N}  U_{{x}_{1}^{+}, \mu}^{2^*-1}  U_{{x}_{j}^{+}, \mu}
= \frac{B_0}{\mu^{N-2}|{x}_{j}^{+}-{x}_{1}^{+}|^{N-2}}+O\Big(\frac 1 {\mu^{N-\sigma}|{x}_{j}^{+}-{x}_{1}^{+}|^{N-\sigma}}\Big),  \quad \text{for} ~ j=2,\cdots, k,
\end{equation}
and
\begin{equation}\label{2*-1+-}
    \int_{\R^N}  U_{{x}_{1}^{+}, \mu}^{2^*-1}  U_{{x}_{j}^{-}, \mu}
=\frac{B_0}{\mu^{N-2}|{x}_{1}^{+}-{x}_{j}^{-}|^{N-2}}+O\Big(\frac 1 {\mu^{N-\sigma}|{x}_{1}^{+}-{x}_{j}^{-}|^{N-\sigma}}\Big),  \quad \text{for} ~ j=1,\cdots, k,
\end{equation}
 where  $B_0= \displaystyle\int_{\mathbb{R}^N}  \frac{C_N^{2^*}}{(1+z^{2})^{\frac{N+2} 2}} $ and  $ \sigma $ is a  constant small enough.

\end{lemma}
\begin{proof}
We denote that  $ \overline{{d}_j} =|{x}_{j}^{+}-{x}_{1}^{+}|, ~ \underline{{d}_j} =|{x}_{1}^{+}-{x}_{j}^{-}| $ for  $j=1, \cdots, k $, and $\overline{B^C_{1j}}= \mathbb{R}^N \setminus \left(B_{\frac{\overline{{d}_j} }{2}}({x}_{1}^{+})\cup B_{\frac{\overline{{d}_j} }{2}}({x}_{j}^{+})\right).$  Then
\begin{equation}\label{a1}
    \begin{aligned}
    &\int_{\mathbb{R}^N}  U_{{x}_{1}^{+}, \mu}^{2^*-1} U_{{x}_{j}^{+}, \mu}\\
 =&\int_{B_{\frac{\overline{{d}_j} }{2}}({x}_{1}^{+})}  U_{{x}_{1}^{+}, \mu}^{2^*-1} U_{{x}_{j}^{+}, \mu}+\int_{B_{\frac{\overline{{d}_j} }{2}}({x}_{j}^{+})}  U_{{x}_{1}^{+}, \mu}^{2^*-1} U_{{x}_{j}^{+}, \mu}+\int_{\overline{B^C_{1j}}}  U_{{x}_{1}^{+}, \mu}^{2^*-1} U_{{x}_{j}^{+}, \mu}\\
 :=&I_1+I_2+I_3.
\end{aligned}
\end{equation}\label{a11}
We first consider $I_1$.
\begin{equation}
    \begin{aligned}
    I_1=&\int_{B_{\frac{\overline{{d}_j} }{2}}({x}_{1}^{+})}  U_{{x}_{1}^{+}, \mu}^{2^*-1} U_{{x}_{j}^{+}, \mu}\\
    =&C_N^{2^*}\int_{B_{\frac{\mu\overline{{d}_j} }{2}}(0)}  \frac{1}{(1+z^{2})^{\frac{N+2} 2}}\\
    &\times\left(  \frac{1}{\mu^{N-2}|{x}_{j}^{+}-{x}_{1}^{+}|^{N-2}}-\frac{N-2}{2}  \frac{1+z^{2}+2 \mu\langle z , {x}_{1}^{+}-{x}_{j}^{+} \rangle}{\mu^{N}|{x}_{j}^{+}-{x}_{1}^{+}|^{N}}+O \Big(\frac{(1+z^{2}+2  \mu\langle z , {x}_{1}^{+}-{x}_{j}^{+} \rangle)^{2}}{\mu^{N+2}|{x}_{j}^{+}-{x}_{1}^{+}|^{N+2}}\Big) \right)\text{d} z\\
    =&\frac{B_0}{\mu^{N-2}|{x}_{j}^{+}-{x}_{1}^{+}|^{N-2}}+O\left(\frac{1}{\mu^{N-\sigma}|{x}_{j}^{+}-{x}_{1}^{+}|^{N-\sigma}}\right).
\end{aligned}
\end{equation}
Then, we calculate $I_2$.
\begin{equation}\label{a12}
    \begin{aligned}
 I_2=&\int_{B_{\frac{\overline{{d}_j} }{2}}({x}_{j}^{+})}  U_{{x}_{1}^{+}, \mu}^{2^*-1} U_{{x}_{j}^{+}, \mu}\\
    =&C_N^{2^*}\int_{B_{\frac{\mu\overline{{d}_j} }{2}}(0)}  \frac{1}{(1+z^{2})^{\frac{N-2} 2}}  \frac{1}{(1+|z+\mu{x}_{j}^{+}-\mu{x}_{1}^{+}|^{2})^{\frac{N+2}{2}}}\text{d} z \\
    =&O\left(\frac{1}{\mu^{N+2}|{x}_{j}^{+}-{x}_{1}^{+}|^{N+2}}\right).
\end{aligned}
\end{equation}
And for $I_3$, we have
\begin{equation}\label{a13}
    \begin{aligned}
 I_3=&\int_{\overline{B^C_{1j}}}  U_{{x}_{1}^{+}, \mu}^{2^*-1} U_{{x}_{j}^{+}, \mu}=O\left(\frac{1}{\mu^{N}|{x}_{j}^{+}-{x}_{1}^{+}|^{N}}\right).
 \end{aligned}
\end{equation}
From \eqref{a1}--\eqref{a13}, we can prove \eqref{2*-1++}. Similarly, we can also prove\eqref{2*-1+-}.
 \end{proof}

\medskip

\begin{lemma}\label{A.2}
   Assume that $\displaystyle\frac{k}{\mu\sqrt{1-h^2}}=o(1), \mu h \to +\infty,N\ge 5$, we have the following expansion, for $k \to +\infty$:
\begin{equation}\label{++}
    \int_{\R^N}  U_{{x}_{1}^{+}, \mu}  U_{{x}_{j}^{+}, \mu}
= O\left(\frac{1}{\mu^{N-2}|{x}_{j}^{+}-{x}_{1}^{+}|^{N-4}}\right), \quad \text{for} ~ j=2,\cdots, k,
\end{equation}
and
\begin{equation}\label{+-}
    \int_{\R^N}  U_{{x}_{1}^{+}, \mu}  U_{{x}_{j}^{-}, \mu}
=O\left(\frac{1}{\mu^{N-2}|{x}_{1}^{+}-{x}_{j}^{-}|^{N-4}}\right),  \quad \text{for} ~ j=1,\cdots, k.
\end{equation}

\end{lemma}
\begin{proof}
We consider that
\begin{equation}\label{a2}
    \begin{aligned}
    &\int_{\mathbb{R}^N}  U_{{x}_{1}^{+}, \mu} U_{{x}_{j}^{+}, \mu}\\
 =&\int_{B_{\frac{\overline{{d}_j} }{2}}({x}_{1}^{+})}  U_{{x}_{1}^{+}, \mu} U_{{x}_{j}^{+}, \mu}+\int_{B_{\frac{\overline{{d}_j} }{2}}({x}_{j}^{+})}  U_{{x}_{1}^{+}, \mu} U_{{x}_{j}^{+}, \mu}+\int_{\overline{B^C_{1j}}}  U_{{x}_{1}^{+}, \mu} U_{{x}_{j}^{+}, \mu}\\
 :=&I_1+I_2+I_3.
\end{aligned}
\end{equation}\label{a21}
We first calculate $I_1$.
\begin{equation}
    \begin{aligned}
    I_1=&\int_{B_{\frac{\overline{{d}_j} }{2}}({x}_{1}^{+})}  U_{{x}_{1}^{+}, \mu} U_{{x}_{j}^{+}, \mu}\\
    =&\frac{C_N^{2}}{\mu^2}\int_{B_{\frac{\mu\overline{{d}_j} }{2}}(0)}  \frac{1}{(1+z^{2})^{\frac{N-2} 2}}  \frac{1}{(1+|z+\mu{x}_{1}^{+}-\mu{x}_{j}^{+}|^{2})^{\frac{N-2}{2}}}\text{d} z \\
    =&O\left(\frac{1}{\mu^{N-2}|{x}_{j}^{+}-{x}_{1}^{+}|^{N-4}}\right).
\end{aligned}
\end{equation}
Similarly, we get
\begin{equation}\label{a22}
    \begin{aligned}
 I_2=O\left(\frac{1}{\mu^{N-2}|{x}_{j}^{+}-{x}_{1}^{+}|^{N-4}}\right).
\end{aligned}
\end{equation}
As for $I_3$, we have
\begin{equation}\label{a23}
    \begin{aligned}
 I_3=&\int_{\overline{B^C_{1j}}}  U_{{x}_{1}^{+}, \mu}U_{{x}_{j}^{+}, \mu}=O\left(\frac{1}{\mu^{N-2}|{x}_{j}^{+}-{x}_{1}^{+}|^{N-4}}\right).
 \end{aligned}
\end{equation}
Combining \eqref{a2}--\eqref{a23}, we get \eqref{++}. Similarly, we can also prove\eqref{+-}.
 \end{proof}
\medskip
\begin{lemma}\label{B.1}
Assume that $\alpha\ge1,N\ge 5$, we have the following estimates for $k \to +\infty, j=2,\cdots,k$:
\begin{equation}\label{b1}
    \begin{aligned}
\sum_{j=2}^k  \frac{1}{|x_{1}^{+}-x_{j}^{+}|^{\alpha}}
& =\left\{\begin{array}{ll}
O\left(\displaystyle\frac{k^\alpha}{\left(r\sqrt{1-h^2}\right)^{\alpha}}\right), & \text { if } \alpha>1 ;\\
\\
O\left(\displaystyle\frac{k^\alpha \ln k}{\left(r\sqrt{1-h^2}\right)^{\alpha}}\right), & \text { if } \alpha=1,
\end{array}\right.
\end{aligned}
\end{equation}
and for $j=1,\cdots,k$, we have:
\begin{equation}\label{b2}
    \begin{aligned}
\sum_{j=1}^k  \frac{1}{|x_{1}^{+}-x_{j}^{-}|^{\alpha}}
&= \left\{\begin{array}{ll}
O\left(\displaystyle\frac{k}{r^\alpha}\right), & \text { if } h \neq o_k(1) ;\\
\\
O\left(\displaystyle\frac{k^{\alpha-1} }{r^{\alpha}h}\right)+O\left(\displaystyle\frac{k}{r^\alpha}\right), & \text { if } hk = o_k(1);\\
\\
O\left(\displaystyle\frac{k^{\alpha }}{r^{\alpha}}\right)=O\left(\displaystyle\frac{1}{r^{\alpha}h^\alpha}\right), & \text { if there exist  }c_2>c_1>0, c_1<hk<c_2;\\
\\
O\left(\displaystyle\frac{k }{r^{\alpha}h^{\alpha-1}}\right), & \text { if } h=o_k(1), hk \to \infty,\alpha>1;\\
\\
O\left(\displaystyle\frac{-k \ln h}{r^{\alpha}}\right), & \text { if } h=o_k(1), hk \to \infty,\alpha=1.
\end{array}\right.
\end{aligned}
\end{equation}

In particular, under the conditions that $\alpha=N-2$, we have the following expansions for $k \to +\infty$:
\begin{align} \label{b3}
\sum_{j=2}^k  \frac{1}{|x_{1}^{+}-x_{j}^{+}|^{N-2}} =  \frac{B_1 k^{N-2}}{\big(r \sqrt{1-h^2}\big)^{N-2}}  \big(1+\zeta_1(k)\big),\quad \text{for} ~ j=2,\cdots, k,
\end{align}
and if $hk \to +\infty$, then
 \begin{equation}
 \begin{split} \label{b4}
& \sum_{j=1}^k  \frac{1}{|x_{1}^{+}-x_{j}^{-}|^{N-2}} = \frac{B_2 k}{r^{N-2} h^{N-3} \sqrt {1-h^2}} \, \big(1+\zeta_2(k)\big)+\frac{\zeta_1(k) k^{N-2}}{\big(r \sqrt{1-h^2}\big)^{N-2}},\quad \text{for} ~ j=1,\cdots, k,
\end{split}
\end{equation}
where
\begin{align} \label{B1B2}
B_1=  \frac 2 {(2\pi)^{N-2}}  \sum_{j=1}^{+\infty} \frac 1 {j^{N-2}},
\quad B_2= \frac 1 {2^{N-3} \pi}\int_0^{+\infty}\,\frac{1}{\big(1+z^2\big)^{\frac{N-2}{2}}}\,{\mathrm d}z,
\end{align}

\begin{align} \label{zeta}
\zeta_1(k)=
\begin{cases}
O\left(\displaystyle\frac{1}{k^{2}}\right),  \ \quad  N\ge 6,
\\
\\
O\left(\displaystyle\frac{\ln k}{k^2}\right), \quad N=5,
\end{cases}  \quad \zeta_2(k)= O\big((hk)^{-1}\big).
\end{align}
\end{lemma}
\begin{proof}
Noting that
\begin{equation}
    |x_{1}^{+}-x_{j}^{+}|=2r\sqrt{1-h^2}\sin\frac{(j-1)\pi}{k},
\end{equation}
and
\begin{equation}
     |x_{1}^{+}-x_{j}^{-}|=2r\left((1-h^2)\left(\sin\frac{(j-1)\pi}{k}\right)^2+h^2\right)^{\frac{1}{2}}.
\end{equation}
Then we have
\begin{equation}\label{++alpha}
    \begin{aligned}
\sum_{j=2}^k  \frac{1}{|x_{1}^{+}-x_{j}^{+}|^{\alpha}}&=\frac{1}{\left(2r\sqrt{1-h^2}\right)^\alpha}\sum_{j=2}^k \frac{1}{\left(\sin\frac{(j-1)\pi}{k}\right)^\alpha} \\
& =\left\{\begin{array}{ll}
\displaystyle\frac{2}{\left(2r\sqrt{1-h^2}\right)^{\alpha}} \displaystyle\sum_{j=2}^{\frac{k}{2}} \frac{1}{\left(\sin \frac{(j-1) \pi}{k}\right)^{\alpha}}+\frac{1}{\left(2r\sqrt{1-h^2}\right)^{\alpha}}, & \text { if } k \text { is even; } \\
\displaystyle\frac{2}{\left(2r\sqrt{1-h^2}\right)^{\alpha}} \displaystyle\sum_{j=2}^{\frac{k+1}{2}} \frac{1}{\left(\sin \frac{(j-1) \pi}{k}\right)^{\alpha}}, & \text { if } k \text { is odd. }
\end{array}\right.
\end{aligned}
\end{equation}
For $j=2,\cdots,[\frac{k+1}{2}]$, we have that
\begin{equation*}
    \frac{2(j-1)}{k}\le \sin \frac{(j-1)\pi}{k}\le \frac{(j-1)\pi}{k},
\end{equation*}
thus, \eqref{b1} holds.

Similarly, we have
\begin{equation}\label{+-alpha}
    \begin{aligned}
\sum_{j=1}^k  \frac{1}{|x_{1}^{+}-x_{j}^{-}|^{\alpha}}
&=\frac{2}{\left(2r\right)^\alpha}\displaystyle\sum_{j=1}^{\left[\frac{k-1}{2}\right]} \frac{1}{\left((1-h^2)\left(\sin\frac{j\pi}{k}\right)^2+h^2\right)^{\frac{\alpha}{2}}}+O\left(\displaystyle\frac{1}{r^\alpha}\right)\\
&= \left\{\begin{array}{ll}
O\left(\displaystyle\frac{k}{r^\alpha}\right), & \text { if } h \neq o_k(1) ;\\
\\
O\left(\displaystyle\frac{k^{\alpha-1} }{r^{\alpha}h}\right)+O\left(\displaystyle\frac{k}{r^\alpha}\right), & \text { if } hk = o_k(1);\\
\\
O\left(\displaystyle\frac{k^{\alpha }}{r^{\alpha}}\right)=O\left(\displaystyle\frac{1}{r^{\alpha}h^\alpha}\right), & \text { if there exist  }c_2>c_1>0, c_1<hk<c_2;\\
\\
O\left(\displaystyle\frac{k }{r^{\alpha}h^{\alpha-1}}\right), & \text { if } h=o_k(1), hk \to \infty,\alpha>1;\\
\\
O\left(\displaystyle\frac{-k \ln h}{r^{\alpha}}\right), & \text { if }h=o_k(1), hk \to \infty,\alpha=1.
\end{array}\right.
\end{aligned}
\end{equation}
The last equation comes from, with the assumption $h=o_k(1)$ and a fixed integer $1<k_1<k$, that

    \begin{align*}
&\displaystyle\sum_{j=1}^{k_1} \frac{1}{\left((1-h^2)\left(\sin\frac{j\pi}{k}\right)^2+h^2\right)^{\frac{\alpha}{2}}}\\&= \left\{\begin{array}{ll}
k^\alpha \displaystyle\sum_{j=1}^{k_1} \frac{1}{\left(1+\frac{h^2}{(1-h^2)\left(\sin\frac{j\pi}{k}\right)^2}\right)^{\frac{\alpha}{2}}}+O(k), & \text { if }  hk = o_k(1);\\
\\
O\left(\displaystyle\frac{1}{h^\alpha}\right)+O(k),  &\text { if there exist } c_2>c_1>0, c_1<hk<c_2;\\
\\
\displaystyle\frac{1}{h^\alpha}\displaystyle\sum_{j=1}^{k_1} \frac{1}{\left(1+\frac{(1-h^2)\left(\sin\frac{j\pi}{k}\right)^2}{h^2}\right)^{\frac{\alpha}{2}}}+O(k), &\text { if }  h=o_k(1),  hk \to \infty;
\end{array}\right.\\
&= \left\{\begin{array}{ll}
k^\alpha O\left(\frac{1}{hk}\displaystyle\int_0^ {+\infty} \frac{1}{(1+z^2)^\frac{\alpha}{2}} \text{d} z\right)+O(k), & \text { if } hk = o_k(1);\\
\\
O\left(\displaystyle\frac{1}{h^\alpha}\right), & \text { if there exist  }c_2>c_1>0, c_1<hk<c_2;\\
\\
\displaystyle\frac{1}{h^\alpha}O\left(hk\int_0^ {+\infty} \frac{1}{(1+z^2)^\frac{\alpha}{2}} \text{d} z\right), & \text { if } h=o_k(1), hk \to \infty, \alpha>1\\
\\
\displaystyle\frac{1}{h}O\left(hk\int_0^ {(1-h^2)\pi^2h^{-2}} \frac{1}{1+z}\text{d} z\right), & \text { if } h=o_k(1), hk \to \infty, \alpha=1.
\end{array}\right.
\end{align*}

Next, we assume that $\alpha=N-2$, then we can calculate
\begin{align*}
   \displaystyle\sum_{j=2}^{\left[\frac{k+1}{2}\right]} \frac{1}{\left(\sin \frac{(j-1) \pi}{k}\right)^{N-2}}&= \displaystyle\sum_{j=2}^{\left[\frac{k}{6}\right]} \frac{1}{\left(\sin \frac{(j-1) \pi}{k}\right)^{N-2}}+ \displaystyle\sum_{j=\left[\frac{k}{6}\right]+1}^{\left[\frac{k+1}{2}\right]} \frac{1}{\left(\sin \frac{(j-1) \pi}{k}\right)^{N-2}} \\
   &=\displaystyle\sum_{j=1}^{\left[\frac{k}{6}\right]-1} \frac{1}{\left(\frac{j \pi}{k}\right)^{{N-2}}}\left(1+O\left(\frac{j^2}{k^2}\right)\right)+O(  k )\\
   &=\frac{k^{N-2}}{\pi^{N-2}}\left({B_{1}'}+\zeta_1(k)\right),
\end{align*}
where  $  {B_{1}'}=   \displaystyle\sum_{j=1}^{+\infty} \frac 1 {j^{N-2}} $
and  $\zeta_1(k)  $ is defined in \eqref{zeta}. Using \eqref{++alpha}, we have proved \eqref{b3}.

On the other hand, assume that $\alpha=N-2, hk \to +\infty$, then
 \begin{equation}
\begin{aligned}
    &\displaystyle\sum_{j=1}^{k} \frac{1}{\left((1-h^2)\left(\sin\frac{(j-1)\pi}{k}\right)^2+h^2\right)^{\frac{N-2}{2}}} \\
=&2\displaystyle\sum_{j=1}^{\left[\frac{k}{6}\right]} \frac{1}{\left((1-h^2)\left(\sin\frac{j\pi}{k}\right)^2+h^2\right)^{\frac{N-2}{2}}}+O\left(\frac{1}{h^{N-2}}\right)+O(k)\\
=&2\displaystyle\sum_{j=1}^{\left[\frac{k-1}{2}\right]} \frac{1}{\left((1-h^2)\left(\frac{j\pi}{k}\right)^2+h^2\right)^{\frac{N-2}{2}}}+O\left(\sum_{j=1}^{\left[\frac{k-1}{2}\right]}\frac{(1-h^2)\left(\frac{j\pi}{k}\right)^4}{\left(\left(1-h^2\right)\left(\frac{j\pi}{k}\right)^2+h^2\right)^\frac{N}{2}}    \right)+O\left(\frac{1}{h^{N-2}}\right)+O(k)\\
=&2\displaystyle\sum_{j=1}^{\left[\frac{k-1}{2}\right]} \frac{1}{\left((1-h^2)\left(\frac{j\pi}{k}\right)^2+h^2\right)^{\frac{N-2}{2}}}+ \zeta_1(k) O\left(  \frac{k^{N-2}}{(\sqrt{1-h^2})^{N-2}}\right)+O\left(\frac{1}{h^{N-2}}\right)+O(k)\\
=&\frac{2k}{h^{N-3}\sqrt{1-h^2}\pi}\int_0^{+\infty} \frac{1}{(1+z^2)^\frac{N-2}{2}} \mathrm{d} z\big(1+O\big((hk)^{-1}\big)+O\left({h^{N-3}}\right)\big)+\zeta_1(k) O\left(  \frac{k^{N-2}}{(\sqrt{1-h^2})^{N-2}}\right).
\end{aligned}
   \end{equation}
Using \eqref{+-alpha}, we can prove \eqref{b4}.
\end{proof}

\begin{lemma}\label{B.2}
Assume that $\alpha\ge1,N\ge 5, (r,h,\mu)\in {\mathscr S}_k,$ we have the following expansions for $k \to +\infty$:
\begin{align} \label{b5}
\sum_{j=2}^k  \frac{\left(1-\cos{\frac{2(j-1)\pi}{k}}\right)}{|x_{1}^{+}-x_{j}^{+}|^{N}} =  \frac{1}{2}\frac{B_1 k^{N-2}}{\left(r\sqrt{1-h^2}\right)^N}\left(1+\zeta_1(k)\right),
\end{align}
 \begin{equation}
 \begin{split} \label{b6}
& \sum_{j=1}^k  \frac{\left(1-\cos{\frac{2(j-1)\pi}{k}}\right)}{|x_{1}^{+}-x_{j}^{-}|^{N}} = \frac{B_2k}{2(N-2)r^N h^{N-3}(\sqrt{1-h^2})^3}+\frac{\zeta_1(k)k^{N-2}}{(\sqrt{1-h^2})^{N-2}},
\end{split}
\end{equation}
and
\begin{equation}
 \begin{split} \label{b7}
& \sum_{j=1}^k  \frac{1}{|x_{1}^{+}-x_{j}^{-}|^{N}} = \frac{(N-3)B_2 k}{(N-2)r^{N} h^{N-1} \sqrt {1-h^2}} \, \big(1+\zeta_2(k)\big)+O\left(\frac{k}{h^{N-1}}\right),
\end{split}
\end{equation}
where
where  $  B_{1}, B_{2} $ are defined in \eqref{B1B2}
and  $\zeta_1(k)  $ is defined in \eqref{zeta}.
\end{lemma}
\begin{proof}
Similar to the calculation in Lemma \ref{B.1}, we have
\begin{equation*}
    \begin{aligned}
\sum_{j=2}^k  \frac{\left(1-\cos{\frac{2(j-1)\pi}{k}}\right)}{|x_{1}^{+}-x_{j}^{+}|^{N}}=\frac{1}{\left(2r\sqrt{1-h^2}\right)^N}\sum_{j=2}^k \frac{\left(1-\cos{\frac{2(j-1)\pi}{k}}\right)}{\left(\sin\frac{(j-1)\pi}{k}\right)^N} 
 = \frac{1}{2}\frac{B_1 k^{N-2}}{\left(r\sqrt{1-h^2}\right)^N}\left(1+\zeta_1(k)\right).
\end{aligned}
\end{equation*}
On the other hand, noting that for $\alpha\ge 1,$
\begin{equation}
    \int_0^{+\infty} \frac{1}{(1+z^2)^\frac{\alpha}{2}} \text{d} z =\frac{\Gamma(\frac{1}{2})\Gamma(\frac{\alpha-1}{2})}{\Gamma(\frac{\alpha}{2})},
    \end{equation}
    where $\Gamma(x)$ is the Gamma function. Then,
\begin{equation*}
    \begin{aligned}
&\sum_{j=1}^k  \frac{\left(1-\cos{\frac{2(j-1)\pi}{k}}\right)}{|x_{1}^{+}-x_{j}^{-}|^{N}}\\
=&\frac{2}{\left(2r\right)^N}\displaystyle\sum_{j=1}^{\left[\frac{k-1}{2}\right]} \frac{\left(1-\cos{\frac{2j\pi}{k}}\right)}{\left((1-h^2)\left(\sin\frac{j\pi}{k}\right)^2+h^2\right)^{\frac{N}{2}}}+O\left(1\right)\\
=&\frac{2}{\left(2r\right)^N} \left(\displaystyle\sum_{j=1}^{\left[\frac{k}{12}\right]} \frac{\left(1-\cos{\frac{2j\pi}{k}}\right)}{\left((1-h^2)\left(\sin\frac{j\pi}{k}\right)^2+h^2\right)^{\frac{N}{2}}}+O\left(k\right)\right)\\
=&\frac{B_2k}{2r^Nh^{N-3}(\sqrt{1-h^2})^3}\left(1-\frac{N-3}{N-2} \right)+O\left(\frac{k}{h^{N-5}}\right)+\zeta_1(k)O\left(\frac{k^{N-2}}{(\sqrt{1-h^2})^{N-2}}\right)\\
=&\frac{B_2k}{2(N-2)r^N h^{N-3}(\sqrt{1-h^2})^3}+\zeta_1(k)O\left(\frac{k^{N-2}}{(\sqrt{1-h^2})^{N-2}}\right).
\end{aligned}
\end{equation*}

Finally, by a similar but more accurate arguments in Lemma \ref{B.1}, we can proof \eqref{b7}.

\end{proof}

\medskip

\begin{lemma}\label{B.3}
    Suppose that $N\ge 5,(r,h,\mu)\in\mathscr{S}_k,$ then for $y\in \Omega^+_1$, ${\alpha_1} \in [1, N-2],$ there exists a constant $C>0$ that
    \begin{equation}
      \left( \sum_{j=2}^k U_{x^{+}_j,\mu}+\sum_{j=1}^kU_{x^{-}_j,\mu}\right) \le \left\{\begin{matrix}
    \displaystyle\frac{C\mu^{\frac{N-2}{2}}}{(1+\mu|y-{x}_{1}^{+}|)^{N-2-{\alpha_1}} }\frac{k^{\alpha_1}}{\mu^{\alpha_1}} \qquad \text{if} \ {\alpha_1}>1, \\
    \\
    \displaystyle\frac{C\mu^{\frac{N-2}{2}}}{(1+\mu|y-{x}_{1}^{+}|)^{N-2-{\alpha_1}} }\frac{k \ln k}{\mu}\qquad \text{if} \ {\alpha_1}=1.
\end{matrix}\right.
    \end{equation}
\end{lemma}
\begin{proof}
    For any $y\in \Omega^+_1$, we have
    \begin{equation*}
        |y-x^{+}_j|\ge \frac{1}{2} |x^{+}_j-x^{+}_1|,\qquad j=2,\cdots,k,
    \end{equation*}
    and
    \begin{equation*}
        |y-x^{-}_j|\ge \frac{1}{2}|x^{+}_1-x^{-}_1|\qquad j=1,\cdots,k.
    \end{equation*}
Then we have
 \begin{equation*}
        \begin{aligned}
\left(\sum_{j=2}^{k}U_{{x}_{j}^{+},\mu}+\sum_{j=1}^{k}U_{{x}_{j}^{-},\mu}\right) \le& \frac{C\mu^{\frac{N-2}{2}}}{(1+\mu^2|y-{x}_{1}^{+}|^2)^\frac{N-2-{\alpha_1}}{2} }\left(\sum_{j=2}^{k}\frac{1}{\mu^{\alpha_1}|{x}_{j}^{+}-{x}_{1}^{+}|^{\alpha_1}}+\sum_{j=1}^{k}\frac{1}{\mu^{\alpha_1}|{x}_{1}^{+}-{x}_{j}^{-}|^{\alpha_1}}\right)\\
\le &\left\{\begin{matrix}
    \displaystyle\frac{C\mu^{\frac{N-2}{2}}}{(1+\mu|y-{x}_{1}^{+}|)^{N-2-{\alpha_1}} }\frac{k^{\alpha_1}}{\mu^{\alpha_1}} \qquad \text{if} \ {\alpha_1}>1, \\
    \\
    \displaystyle\frac{C\mu^{\frac{N-2}{2}}}{(1+\mu|y-{x}_{1}^{+}|)^{N-2-{\alpha_1}} }\frac{k \ln k}{\mu}\qquad \text{if} \ {\alpha_1}=1.
\end{matrix}\right.
        \end{aligned}
    \end{equation*}

\end{proof}

\medskip

\begin{lemma}\label{B.31}
    Suppose that $N\ge 5,(r,h,\mu)\in\mathscr{S}_k,$ then for $y\in \Omega^+_1$, ${\alpha_2} \in [1, N-1],$ there exists a constant $C>0$ that
    \begin{equation}
      \left( \sum_{j=2}^k \overline{\mathbb Z}_{2j}+\sum_{j=1}^k\underline{\mathbb Z}_{2j}\right) \le \left\{\begin{matrix}
    \displaystyle\frac{C\mu^{\frac{N}{2}}}{(1+\mu|y-{x}_{1}^{+}|)^{N-1-{\alpha_2}} }\frac{k^{\alpha_2}}{\mu^{\alpha_2}} \qquad \text{if} \ {\alpha_2}>1, \\
    \\
    \displaystyle\frac{C\mu^{\frac{N}{2}}}{(1+\mu|y-{x}_{1}^{+}|)^{N-1-{\alpha_2}} }\frac{k \ln k}{\mu}\qquad \text{if} \ {\alpha_2}=1.
\end{matrix}\right.
    \end{equation}
\end{lemma}
\begin{proof}
    The proof is similar to that of Lemma \ref{B.3}. Here we omit it.
    
\end{proof}

\end{document}